\documentclass[11pt]{article}

\usepackage[colorlinks=true, linkcolor=black, anchorcolor= black, citecolor= black, filecolor= black, menucolor= black, pagecolor= black, urlcolor=red]{hyperref}

\newtheorem{thm}{Theorem}[subsection]
\newtheorem{prop}[thm]{Proposition}
\newtheorem{lem}[thm]{Lemma}
\newtheorem{cor}[thm]{Corollary}

\newtheorem{defn}[thm]{Definition}

\newtheorem{example}[thm]{Example}
\newtheorem{remark}[thm]{Remark}
\newtheorem{conj}[thm]{Conjecture}
\newtheorem{prob}[thm]{Problem}

\newcommand{\pr}[1]{{{\bf P}^{#1}}}
\newcommand{\skipit}[1]{{}}
\newcommand{\prfend}{\hbox to7pt{\hfil}
\par\vskip-\baselineskip\hbox to\hsize
{\hfil\vbox {\hrule width6pt height6pt}}\vskip\baselineskip}

\newcommand{\myarrow}[2]{\hbox to #1pt{\hfil$\to$\hfil}{\hskip-#1pt{\raise
10pt\hbox to#1pt{\hfil$\scriptscriptstyle #2$\hfil}}}}

\long\def\eatit#1{}

\begin{document}
\title{Betti numbers for fat point ideals in the plane: a geometric approach.\thanks{{\it 2000 Mathematics Subject Classification}. Primary 14C20, 13P10; Secondary 14J26, 14J60.\newline
\hbox to.2in{\hfil} {\it Key words and phrases}. Graded Betti numbers, fat points, splitting types.}}

\author{Alessandro Gimigliano\\
Dipartimento di Matematica e CIRAM\\
Universit\`a di Bologna\\
40126 Bologna, Italy\\
email: gimiglia@dm.unibo.it
\and
Brian Harbourne\\
Department of Mathematics\\
University of Nebraska\\
Lincoln, NE 68588-0130 USA\\
email: bharbour@math.unl.edu
\and
Monica Id\`a\\
Dipartimento di Matematica\\
Universit\`a di Bologna\\
40126 Bologna, Italy\\
email: ida@dm.unibo.it}

\date{June 15, 2007}

\maketitle

\thanks{Acknowledgments: We thank GNSAGA, and
the University of Bologna, which supported
visits to Bologna by the second author, who also
thanks the NSA and NSF for supporting his research.
We also thank the referee for his careful and helpful 
comments.}

\begin{abstract}
We consider the open problem of determining the
graded Betti numbers for fat point subschemes $Z$
supported at general points of $\pr2$.
We relate this problem to the open geometric problem
of determining the splitting type of the pullback
of $\Omega_{\pr2}$ to the normalization of certain 
rational plane curves. We give a conjecture 
for the graded Betti numbers which 
would determine them in all degrees but 
one for every fat point subscheme
supported at general points of $\pr2$.
We also prove our Betti number
conjecture in a broad range of cases. An appendix 
discusses many more cases in which
our conjecture has been verified computationally
and provides a new and more efficient computational
approach for computing graded Betti numbers
in certain degrees. It also demonstrates how to derive explicit
conjectural values for the Betti numbers and 
how to compute splitting types.
\end{abstract}

\section{Introduction}\label{intro}

Let $Z=m_1P_1+\cdots+m_nP_n$ be a fat point
subscheme of $\pr2$ supported at general points $P_i$.
Thus $Z$ is the 0-dimensional subscheme of $\pr2$
defined by the homogeneous ideal $I(Z)=\cap_i I(P_i)^{m_i}$ 
in the homogeneous coordinate ring $R=K[\pr2]$ of $\pr2$
(where we take $K$ to be an algebraically closed field of arbitrary 
characteristic), where $I(P_i)$ is the ideal 
generated by all homogeneous forms $f\in R$ vanishing at $P_i$. 
The homogeneous component $I(Z)_t$ of $I(Z)$ in degree $t$
is just the $K$-vector space span of the homogeneous elements of
$I(Z)$ of degree $t$. Thus $I(Z)_t$ consists of all 
homogeneous polynomials of degree $t$ which vanish at each point $P_i$
to order at least $m_i$. The Hilbert function of $I(Z)$
is the function $h_Z$ which gives the $K$-vector space dimension  
$h_Z(t)=\hbox{dim }I(Z)_t$ of $I(Z)_t$ as a function of $t$.

\subsection{The SHGH Conjecture}

The problem of determining $h_Z$ has attracted a lot of attention over 
the years, and it is still an
open problem in general. In fact, even the most fundamental
problem is open: it is not known in general
what the least $t$ is (which we call $\alpha(Z)$) for which $I(Z)_t\ne0$.
The fact that even the latter problem is open is less surprising
given that if one knows $\alpha(Z)$ for every $Z$, then one can determine
$h_Z(t)$ for all $t$ and $Z$, and conversely.

Given that it is not known how big $I(Z)_t$ is, or even if it is
non-zero, it is not surprising that the least value of $t$
such that $I(Z)$ is generated in degrees $t$ or less is not known. More precisely,
the graded Betti numbers for the minimal free resolution of $I(Z)$
are not known. There is not even a conjecture in general for the Betti numbers.
On the other hand, there is a general conjecture for the values of
the Hilbert function $h_Z$. This is the SHGH Conjecture 
(due, in various equivalent forms,
to Segre \cite{refSe}, Harbourne \cite{refVanc}, 
Gimigliano \cite{refGi} and Hirschowitz \cite{refHi}; see
Conjecture \ref*{SHGHconjVanc} or Conjecture \ref*{SHGHconj}).
Considerations of geometry lead to a lower bound $e(h_Z,t)$ for $h_Z(t)$
(the definition of $e(h_Z,t)$, which is somewhat complicated,
is given in the appendix). The SHGH Conjecture
asserts that $e(h_Z,t)=h_Z(t)$. Thus $e(h_Z,t)$ is 
regarded as the ``expected'' value of $h_Z(t)$.

In degrees $t>\alpha(Z)$, the SHGH Conjecture
implies the following simple statement:

\begin{conj}\label{simpleSHGH1}
Let $Z=m_1P_1+\cdots+m_nP_n$,
for nonnegative integers $m_i$ and general points $P_i\in \pr2$.
If $t>\alpha(Z)$ and $f_1^{b_1}\cdots f_r^{b_r}$
is a factorization of the greatest common divisor
of $I(Z)_t$ as a product of non-associate irreducible factors, then 
$$h_Z(t)={t+2\choose 2}-\sum_i{m_i+1\choose 2}+\sum_l{b_l\choose 2}.$$
\end{conj}

In particular, under the assumptions given,
$e(h_Z,t)={t+2\choose 2}-\sum_i{m_i+1\choose 2}+\sum_l{b_l\choose 2}$.
Although it is perhaps not clear from this how to actually compute $e(h_Z,t)$, 
the full SHGH Conjecture
allows one to compute $e(h_Z,t)$ in terms only of the
$m_i$, and also to conjecturally determine
the degrees and multiplicities at each point $P_i$ of
all curves in the base locus of $I(Z)_t$. \eatit{Refer to where this is discussed; Appendix?}
The advantage of the statement above is that it gives insight 
into the problem of determining $h_Z$, without the burden
of the technicalities needed to state the full SHGH Conjecture.

For example, $I(Z)_t$ consists of all homogeneous polynomials of degree $t$
which vanish at each point $P_i$ to order at least $m_i$.
But the vector space of all forms of degree $t$ has dimension 
${t+2\choose 2}$, and requiring vanishing at $P_i$ to order $m_i$ imposes
${m_i+1\choose 2}$ independent linear conditions. We do not know
that the conditions imposed at one point are independent of 
those imposed at all of the other points,
and in fact they are not always independent. Thus we obtain the bound
$h_Z(t)\ge{t+2\choose 2}-\sum_i{m_i+1\choose 2}$. What the SHGH Conjecture 
does in essence is to give a precise measure of 
the failure of the imposed conditions to be independent.
The key insight is that the failure of independence in degree $t>\alpha(Z)$
is due to the gcd not being square-free. What the conjecture above
hides is that the curves defined by the forms $f_i$ are 
thought always to be very special, and this is significant 
for what we do in this paper.

To make this clearer we strengthen the conjecture above to include
the case that $t=\alpha(Z)$. To do this we need a minor technicality.
We say a plane curve 
$C$ defined by an irreducible form of degree $d$ is 
{\it contributory} with respect to points $P_1,\cdots,P_n$
if $C$ is a rational curve smooth except possibly at the points $P_i$,
such that $\hbox{mult}_{P_1}(C)+\cdots+\hbox{mult}_{P_n}(C)=3d-1$.
We will say a plane curve 
$C$ defined by an irreducible form of degree $d$ is 
{\it negative} with respect to points $P_1,\cdots,P_n$ if
$(\hbox{mult}_{P_1}(C))^2+\cdots+(\hbox{mult}_{P_n}(C))^2>d^2$.
It follows by the genus formula that curves contributory
for $P_1,\ldots,P_n$ are also negative.

After extending Conjecture \ref*{simpleSHGH1}, we have:

\begin{conj}\label{simpleSHGH2}
Let $Z=m_1P_1+\cdots+m_nP_n$,
for nonnegative integers $m_i$ and general points $P_i\in\pr2$.
Given $t\ge\alpha(Z)$, let $f_1^{c_1}\cdots f_r^{c_r}$
be a factorization of the greatest common divisor
of $I(Z)_t$ as a product of non-associate irreducible factors.
Then every factor $f_j$ which 
defines a curve negative for the points $P_i$ defines a 
contributory curve, and we have
$$h_Z(t)={t+2\choose 2}-\sum_i{m_i+1\choose 2}+\sum_k{c_{j_k}\choose 2},$$
where the index $k$ runs over the factors $f_{j_k}$ defining
curves contributory for the points $P_i$.
\end{conj}

The point is that the SHGH Conjecture implies 
that it is only exponents $c_j$ of factors $f_j$ defining
contributory curves which contribute to the lack of independence.
The SHGH Conjecture also implies in degrees $t>\alpha(Z)$ 
that every factor $f_j$ defines a contributory curve, and that
every curve negative for general points is in fact contributory.
(Factors defining non-contributory curves do occur in degree $\alpha$, however;
for example, if $Z=P_1+\cdots+P_9$ consists of 9 general points, 
then $\alpha(Z)=3$ and the gcd in degree 3 is an irreducible cubic
defining an elliptic curve $C$, which is therefore not contributory,
but it is also not negative.)

It is known, however, that contributory curves do contribute to the failure
of independence.
The SHGH Conjecture is that nothing else contributes.
In particular, if the greatest common divisor of $I(Z)_t$ is 1
(i.e., if the zero locus of $I(Z)_t$ is at most 0-dimensional),
or more generally if the gcd of $I(Z)_t$ is just square-free,
then the conjecture is that there is no failure of
independence, and hence that $h_Z(t)={t+2\choose 2}-\sum_i{m_i+1\choose 2}$.

\subsection{The Betti Number Conjecture}\label{bnc}

The goal of this paper is to give a theorem and a conjecture
for graded Betti numbers for ideals $I(Z)$, mimicking 
the fact that there is a bound $e(h_Z,t)\le h_Z(t)$ which by the SHGH 
Conjecture is an equality, given explicitly by Conjecture \ref*{simpleSHGH2}.

The SHGH Conjecture specifies how big the ideal $I(Z)$ is in each degree.
The next question, for which no general conjecture has yet been posed,
is where must one look for generators of $I(Z)$. More precisely,
in any minimal set of homogeneous generators, how many generators are 
there in each degree? The goal of this paper is to develop
(and prove cases of) a conjecture for the numbers of generators
in each degree bigger than $\alpha(Z)+1$. The number of generators in degrees
less than $\alpha(Z)+1$ is trivial: there are obviously no generators in degrees $t<\alpha(Z)$,
and, since any minimal set of 
homogeneous generators of $I(Z)$ must include a $K$-vector 
space basis of $I(Z)_\alpha$,
precisely $h_Z(t)$ generators in degree $t=\alpha(Z)$.
There remains the question of how many generators there are in degree 
$t=\alpha(Z)+1$. Our approach is to relate the number of generators
in a given degree to the splitting of 
a certain rank 2 bundle on certain curves. In degree
$\alpha(Z)+1$ precisely what curves must be taken into account is 
more subtle than it is in larger degrees. Thus here we focus on
degrees larger than $\alpha(Z)+1$. We study the subtleties
needed for a unified approach that subsumes degree $\alpha(Z)+1$ in separate papers,
beginning with \cite{refBMS1}.

For any $t$, the number of homogeneous generators in degree $t+1$
in any minimal set of homogeneous generators is the dimension
of the cokernel of the map $\mu_t:I(Z)_t\otimes R_1\to I(Z)_{t+1}$
given on simple tensors by multiplication, $f\otimes g\mapsto fg$.
What we will do is to give a conjecture for $\hbox{dim cok }\mu_t$
for all $Z$ in all degrees $t>\alpha(Z)$. Clearly
$\hbox{dim cok }\mu_t=0$ if $t<\alpha(Z)-1$ and
$\hbox{dim cok }\mu_t=h_Z(\alpha(Z))$ if $t=\alpha(Z)-1$, so 
our conjecture handles all cases except $t=\alpha(Z)$.
The conjecture is in terms of data determined by contributory
curves.

Let $C'$ be a curve contributory for points $P_1,\ldots,P_n$.
The composition of the normalization map  $C\to C'$
with the inclusion $C'\subset\pr2$ gives a morphism $f:C\to \pr2$.
Pulling back the twisted cotangent bundle $\Omega_{\pr2}(1)$ gives a rank two
bundle $f^*(\Omega_{\pr2}(1))$ on $C$, but $C$ is smooth and rational,
so $f^*(\Omega_{\pr2}(1))$ splits as $f^*(\Omega_{\pr2}(1))\cong 
{\cal O}_{C}(-a_{C})\oplus {\cal O}_{C}(-b_{C})$ for some integers
$a_C\le b_C$. We call $(a_{C}, b_{C})$ the {\it splitting type\/} of $C$
(and for convenience, we set $a_{C'}=a_{C}$ and $b_{C'}=b_{C}$
and refer to $(a_{C'}, b_{C'})$ also as the {\it splitting type\/} of $C'$).

Given $t>\alpha(Z)$, clearly $I(Z)_{t-1}\ne0$, 
so $I(Z)_{t-1}$ has a well-defined gcd (up to scalar multiple).
Let $\gamma_{t-1}$ be the product of all irreducible factors of the gcd
defining curves negative for the points $P_1,\ldots,P_n$.
(If $t>\alpha(Z)+1$, then, as in Conjecture \ref*{simpleSHGH1},
the SHGH Conjecture implies that $\gamma_{t-1}$ 
is itself the gcd.)
Let $\gamma_{t-1}=f_1^{c_1}\cdots f_r^{c_r}$ be its
factorization into non-associate irreducible factors 
and let $d_j=\hbox{deg}(f_j)$.
Do likewise for $I(Z)_t$; it is easy to see that
we get $\gamma_t=f_1^{c_1'}\cdots f_r^{c_r'}$
where $c_j'\le c_j$ for all $j$.
Doing the same for $I(Z)_{t+1}$ gives
$\gamma_{t+1}=f_1^{c_1''}\cdots f_r^{c_r''}$ with $c_j''\le c_j'$.
Let $C_j$ be the normalization of the curve defined by $f_j$.
The following theorem refines ideas of Fitchett \cite{refF1}, \cite{refF2}.

\begin{thm}\label{Bettibound}
Let $Z=m_1P_1+\cdots+m_nP_n$,
for nonnegative integers $m_i$ and general points $P_i\in\pr2$.
Let $t>\alpha(Z)$ with $c_j$, $c_j'$ and $c_j''$ defined as above.
If the SHGH Conjecture holds, then 
$$\hbox{dim cok }\mu_t\le \sum_jd_j(c_j'-c_j'')-\sum_j{c_j'-c_j''\choose 2} 
+ \sum_j\Bigg({c_j-c_j'-a_{C_j}\choose 2}
+{c_j-c_j'-b_{C_j}\choose 2}\Bigg).$$
\end{thm}

In fact, the full SHGH Conjecture is not required
for this theorem. One just needs that $I(Z)$ behaves as expected
for the specific $Z$ being considered.
We also propose the following conjecture:

\begin{conj}\label{Betticonj}
Equality holds in Theorem \ref*{Bettibound}.
\end{conj}

As we show by examples in the appendix, the 
SHGH Conjecture allows one to determine conjectural values for
the exponents $c_j$,
$c_j'$ and $c_j''$. As we show below, in many cases the
splitting type of the curves $C_j$ are also known, and even in those cases where
the type is not known, it is much easier to compute the splitting type
of each $C_j$ symbolically than it is to compute  $\hbox{dim cok }\mu_t$
symbolically in the usual way of finding a
Gr\"obner basis of the ideal $I(Z)$. Thus 
the SHGH Conjecture and Conjecture \ref*{Betticonj}, if true, allow one
to determine the minimal number of homogeneous generators of $I(Z)$
in every degree except possibly degree $\alpha(Z)+1$. Examples in 
the appendix show how this is done. 

We now put this into the context of minimal free graded resolutions.
The minimal free graded resolution of $I(Z)$ is an exact sequence
of the form $0\to M_1\to M_0\to I(Z)\to0$, where
$M_0$ (the module of generators)
and $M_1$ (the module of syzygies) 
are free graded $R$-modules, hence of the form
$M_0=\oplus_{i\ge0} R[-i]^{g_i(Z)}$, and 
$M_1=\oplus_{j\ge0} R[-j]^{s_j(Z)}$
for nonnegative integers $g_i(Z)$ and $s_j(Z)$. (By $R[-i]$
we just mean the free $R$-module of rank 1 with the grading
such that $R[-i]_t=R_{t-i}$.) The graded Betti numbers of $I(Z)$ are
the sequences of integers $g_i(Z)$ and $s_j(Z)$ (which we write as
$g_i$ and $s_j$ if $Z$ is understood). 
The Betti number $g_i$ is just $\hbox{dim cok }\mu_{i-1}$,
hence $g_i$ is the number of generators of degree $i$
in any minimal set of homogeneous generators of $I(Z)$.
Moreover, it is not hard to show that
$g_i-s_i=\Delta^3h_Z(i)$ for all $i$, where $\Delta$ 
is the difference operator; i.e., $\Delta h_Z(i)=h_Z(i)-h_Z(i-1)$
(see p. 685 of \cite{refFHH}).

\subsection{The Structure of The Paper}

We obtain our results by reformulating them in terms of
complete linear systems on the surface $X$ obtained
by blowing up $\pr2$ at the points $P_1,\ldots,P_n$.

In Section \ref*{bckgnd} we recall the background necessary for this 
reformulation, and we state known results needed for our approach.
In Section \ref*{mainresults} we state our main results and show how they
lead to the statement in terms of fat points given above.

We also include an appendix for the purpose of showing
how to obtain explicit predictions for the 
values of Hilbert functions of fat points
and how our results lead to explicit numerical predictions for 
Betti numbers. 
In addition, Section \ref*{algrthms} of the appendix
discusses how to compute splitting types and  
discusses evidence in support of Conjecture \ref*{Betticonj},
partly based on an approach for computing
$g_{i}$ for $i>\alpha+1$ which is substantially more
efficient than the usual methods, which involve finding a Gr\"obner basis
of $I(Z)$. (A Macaulay 2 script which implements
our method is included in Section A2.3
of the posted version of the paper, \cite{refGHIpv}.)
Even for relatively small values of the multiplicities
$m_i$ and even for randomly chosen 
points $P_i$ over a finite field rather than for
generic $P_i$, finding the graded Betti numbers $g_i$ for
$I(Z)$ is beyond what can be done computationally by the usual methods.
For example, we were unable to determine the graded Betti numbers of $I(Z)$
for $Z=77P_1+ \cdots +77P_7 +44P_8+ 11P_9+ 11P_{10}+ 11P_{11}$
by the usual methods, but our new computational method,
based on the results of Section \ref*{mainresults},
working on an 800MHz computer over the finite field $|K|=31991$ 
using randomly chosen points $P_i$,
determined the result in slightly over 5 minutes. The result, of course,
is in agreement with our conjectural expected values.
See also Example \ref*{expex1}.

\subsection{What was Previously Known}

A lot of work has been done on the SHGH Conjecture.
That the SHGH Conjecture holds for $n\le 9$ points was known
to Castelnuovo \cite{refCast}; a more modern proof is given by
Nagata \cite{refN2}. The uniform case
(i.e., $Z=m(P_1+\cdots+P_n)$) was proved 
for $m=2$ by \cite{refAC}, $m=3$ by \cite{refHi2}, for $m\le12$
by \cite{refCM2} and for $m\le 20$ by \cite{refCCMO}. 
The case that $n$ is a square has likewise
seen progressive improvements, with the main difficulty
being to show $I(Z)_t=0$ when it is expected to be.
For example, by specializing $n=16$ points to a 
smooth curve of degree $4$, it is not hard to show that
$h_Z(t)$ has its expected value of ${t+2\choose 2}-4{m+1\choose 2}$ for all 
$t\ge 4m+1$, while $h_Z(t)$ has its expected value
of 0 for $t<4m+1$ by \cite{refN1} for all $m$.
A generalization of this in \cite{refHHF} shows that the SHGH Conjecture holds
in all degrees $t$ such that $e(h_z,t)>0$,
if $n$ is any square as long as $m$ is not too small.
By more technical arguments one can show that
the SHGH Conjecture holds also for small $m$ and $t$.
For example, if $n$ is a power of 4,
\cite{refEvain1} showed SHGH holds in the uniform case;
\cite{refBZ} extended this to $n$ being a
product of powers of 4 and 9; and \cite{refHR} showed
that the SHGH Conjecture holds
for infinitely many $m$ for each square $n$.
By \cite{refEvain2} it is now known to hold in the uniform case
for any $m$ when $n$ is a square; alternate proofs have been given by
\cite{refCM3} and \cite{refR}.

Results for the general case (i.e., such that the multiplicities $m_i$
of $Z=\sum_im_iP_i$ need not all be equal) are not as comprehensive.
That the SHGH Conjecture holds in the case that $m_i\le3$ for all $i$ is 
due to  \cite{refCM1}, improved to $m_i\le4$ by \cite{refmig} and then
to $m_i\le 7$ by \cite{refyang}.

Previous results on graded Betti numbers seem to
start with \cite{refCat}, which obtained a complete answer for $n\le 5$
general points.
This was extended to $n=6$ by \cite{refF3}, then 7 by \cite{refCJM} and 8 by
\cite{refFHH}. For $n>8$, almost all results (and even conjectures)
are for cases which are either uniform or close to uniform.
For example, conjectures in the uniform case were put forward by
\cite{refigp}, and in cases close to uniform by
\cite{refHHF}. (Those conjectures are consistent with 
Conjecture \ref*{Betticonj} due to Corollary \ref*{CMlem}.)
The Betti numbers in the case of $n$ general points of multiplicity $m=1$
were determined by \cite{refGGR}. The uniform Betti numbers conjecture of
\cite{refigp} was verified for $m=2$ by \cite{refI} and for $m=3$ by \cite{refGI}. More generally, if $Z=m_1P_1+\cdots+m_nP_n$
where the points $P_i\in\pr2$ are general and $m_i\le3$ for all $i$,
\cite{refBI} determines the graded Betti numbers in all degrees. 
By Theorem 3.2 of \cite{refHHF}, applying \cite{refEvain2},
it follows that the uniform Betti numbers  conjecture of 
\cite{refigp} holds for all $m\ge (\sqrt{n}-2)/4$ when the number $n$ of points
is an even square. Additional cases are shown in \cite{refHR} when $n$ is 
not a square.

\section{Background}\label{bckgnd}

In this section we set notation and cite well known facts
which we will refer to later in the paper.
Let $P_1,\ldots,P_n$ be distinct (not necessarily general) points of the
projective plane $\pr2$.
Let $p: X\to \pr2$ be the birational morphism given by blowing 
up the points. 

\subsection{Preliminaries}\label{prelims}

The divisor class group $\hbox{Cl}(X)$ of
divisors on $X$ modulo linear equivalence is the free abelian
group with basis $L, E_1,\ldots,E_n$, where $E_i$ is the 
class of the divisor $p^{-1}(P_i)$ and $L$ is the pullback of the class of a line.
Given any divisor $F$ on $X$, the dimension $h^0(X, {\cal O}_X(F))$
of the global sections of ${\cal O}_X(F)$ depends only on
the class $[F]$ of $F$. For convenience, we will denote
$h^0(X, {\cal O}_X(F))$ by either $h^0(X, F)$ or $h^0(X, [F])$,
or even $h^0(F)$ or $h^0([F])$ if $X$ is understood.

Given any $F=tL-m_1E_1-\cdots-m_nE_n$,
by Riemann-Roch we have
$$h^0(X, F)-h^1(X, F)+h^2(F,X)={F^2-K_X\cdot F\over 2}+1,$$
where $K_X=-3L+E_1+\cdots+E_n$ is the canonical class.
Since $E_i$ is reduced and irreducible and $m_i=E_i\cdot F$,
we have a canonical isomorphism $H^0(X, tL-\sum_{m_i>0}m_iE_i)\to H^0(X, F)$.
(The idea is that if $|F|$ is nonempty, then 
$-\sum_{m_i<0}m_iE_i$ is contained in the base locus
of $|F|$, essentially by Bezout's Theorem,
and so $|F|=(-\sum_{m_i<0}m_iE_i)+|tL-\sum_{m_i>0}m_iE_i|$.)
On the other hand, $L$ is nef (meaning that $L\cdot C\ge0$ for any effective
divisor $C$ on $X$), so $h^0(X, F)=0$ if $t<0$.
By duality we have $h^2(X, F)=h^0(X, K_X-F)$, so
it follows that $h^2(X, F)=0$ whenever $t\ge 0$
(in fact, whenever $t\ge -2$). 

If $t\ge0$ and $m_i\ge0$ for all $i$,
then ${F^2-K_X\cdot F\over 2}+1={t+2\choose 2}-
\sum_i {m_i+1\choose 2}$, so Riemann-Roch gives
$$h^0(X, F)\ge \hbox{max}\Bigl(0, {\textstyle{t+2\choose 2}}-
\sum_i {\textstyle {m_i+1\choose 2}}\Bigr).$$
This is just a manifestation of the canonical identification
$H^0(X, tL-m_1E_1-\cdots-m_nE_n)=H^0(\pr2,{\cal I}_Z(t))=I(Z)_t$, where
$Z=m_1P_1+\cdots+m_nP_n$ and ${\cal I}_Z$ is the sheaf of ideals
defining $Z$. Because of this, given any integer $t$ and a fat point scheme
$Z=m_1P_1+\cdots+m_nP_n$, we define
$F_t(Z)=tL-m_1E_1-\cdots-m_nE_n$, and hence
we have $h_Z(t)=h^0(X, F_t(Z))$ for all $t$.

Given a divisor $F$ on $X$, let
$\mu_F:H^0(X, F)\otimes H^0(X, L)\to H^0(X, F+L)$
denote the obvious natural map. By identifying
$H^0(X, F_t(Z))$ with $I(Z)_t$ and $H^0(X, L)$ with 
$H^0(\pr2,{\cal O}_{\pr2}(1))=R_1$, 
it follows that $\hbox{dim}\,\hbox{cok}\,\mu_t(Z)=
\hbox{dim}\,\hbox{cok}\,\mu_{F_t(Z)}$.
Given a sum $F=H+N$ of effective divisors such that $|F|=N+|H|$
(i.e., such that $N$ is contained in the scheme theoretic
base locus of $F$), Lemma \ref*{baselocuslem} gives 
a simple but useful fact which relates 
$\hbox{dim}\,\hbox{cok}\,\mu_{F}$ to $\hbox{dim}\,\hbox{cok}\,\mu_{H}$.
(For the proof, observe that $N$ is in the base locus of 
the image $\hbox{Im}\,\mu_{F}$ of $\mu_{F}$;
i.e., that $\hbox{Im}\,\mu_{F}=N+\hbox{Im}\,\mu_{H}$.
See Lemma 2.10(b) of \cite{refFreeRes}.)
Obviously $\mu_F$ is injective when $F$ is not
effective, and, when $F$ is effective, $F$  decomposes as
in Lemma \ref*{baselocuslem}. Thus
Lemma \ref*{baselocuslem} reduces the general problem of computing
$\hbox{dim}\,\hbox{cok}\,\mu_{F}$ to the case that $|F|$ is effective
and fixed component free, and thus in particular to the case that $F$ is nef.

\begin{lem}\label{baselocuslem}
Let $F=H+N$ be a sum of effective divisors $H$ and $N$ on 
the surface $X$ such that $|F|=N+|H|$. Then
$$\hbox{dim}\,\hbox{cok}\,\mu_{F}=(\hbox{dim}\,\hbox{cok}\,\mu_{H})+(h^0(X,F+L)-h^0(X,H+L)).$$
\end{lem}

It is easy to give examples such that $F=N$ (and hence $\mu_F$ is injective)
but that $N+|L|$ is a proper subset of $|F+L|$ (and hence 
the map $\mu_F$ cannot be surjective). 
For example, if $C$ is an {\it exceptional\/} curve (i.e.,
$C$ is smooth and rational with $C^2=-1$) with $d=C\cdot L$,
then the kernel of $\mu_{L+(d+1)C}$ is non-zero since already $\mu_L$
has non-zero kernel, and $\mu_{L+(d+1)C}$ is not onto since
$C$ is in the base locus of $|L+(d+1)C|$ but $|2L+(d+1)C|$ is base curve free.
In particular, the occurrence of fixed 
components is one reason that $\mu_F$ can fail to have maximal rank
(i.e., fail to be either injective or surjective). 
What motivates this paper is that this is not the only reason.
In fact, $\mu_F$ can fail to have maximal rank even when $F$ is very ample.
(For example, let $F=(3L-E_1-\cdots-E_7)+m(8L-3E_1-\cdots-3E_7-E_8)$,
where the points $P_i$ are general and $m\ge1$. Then $F$ is very ample
\cite{refMA}, but $\mu_F$ fails to have maximal rank \cite{refFHH}.)
The point of this paper is that, when $F-L$ is effective, 
the failure of $\mu_F$ to have maximal rank
depends on the fixed components of $|F-L|$.

Note that if $C$ is a plane curve contributory for points $P_i$,
then $C'$ is an exceptional curve, where $C'$ is
the proper transform of $C$ on the surface $X$ 
obtained by blowing up the points $P_i$.
As discussed in the introduction, our geometric approach for 
determining the dimension of the cokernel of
$\mu_F$ for certain divisors $F$ thus depends on 
knowing the splitting ${\cal O}_E(-a_E)\oplus{\cal O}_E(-b_E)$
of the restriction $p^*\Omega(1)|_E$ of $p^*\Omega(1)$ to exceptional curves $E$.
The following result (see \cite{refAs}, or \cite{refF1}, \cite{refF2})
covers most of what is known:

\begin{lem}\label{splitlem} Let $E\subset X$ be a smooth rational curve,
where $p:X\to \pr2$ is the morphism
blowing up distinct points $P_i$ of $\pr2$. 
Let $d=E\cdot L$ and let $m$ be the
maximum of $E\cdot E_i$, $1\le i\le n$.
Then there are integers
$0\le a_E\le b_E\le d$ with $\hbox{min}(m,d-m)\le a_E\le d-m$ and $d=a_E+b_E$
such that $(p^*\Omega_{\pr2}(1))|_E$
is isomorphic to ${\cal O}_E(-a_E)\oplus{\cal O}_E(-b_E)$.
\end{lem}

Note that, if $d\leq 2m+1$, then $a_E=\hbox{min}(m,d-m)$ and
$b_E=\hbox{max}(m,d-m)$.
For cases not covered by Lemma \ref*{splitlem}, $a_E$ and $b_E$ can be
computed fairly efficiently. We will describe an algorithm for doing so
in section \ref*{algrthms} of the appendix.

\subsection{The SHGH Conjecture}\label{appSHGHsect}

Here we state the version of the SHGH Conjecture given
in \cite{refVanc}. This version is simple to state and useful conceptually.
We include in the appendix an equivalent
version that is more useful for obtaining explicit 
conjectural values of Hilbert functions.

\begin{conj}\label{SHGHconjVanc}
Let $X$ be a surface obtained by blowing up $n$ generic points of $\pr2$.
Then every reduced irreducible curve $C\subset X$ with
$C^2<0$ is an exceptional curve and either $h^0(X, F)=0$ or $h^1(X,F)=0$
for every nef divisor $F$ on $X$.
\end{conj}

\begin{thm}\label{conjequiv}
Conjecture \ref*{SHGHconjVanc} implies 
Conjectures \ref*{simpleSHGH1} and \ref*{simpleSHGH2}.
\end{thm}

\noindent{\bf Proof.} 
First, note that if $H$ is effective, then $h^2(X,H)=0$ by duality
(since $L$ is nef but $(K_X-H)\cdot L<0$). Now let $C$ be nef and effective.
Then $H+C$ is also nef and effective. By Conjecture \ref*{SHGHconjVanc},
$h^1(X,H)=0=h^1(X,H+C)$ and $h^1(X,C)=0=h^2(X,C)$ 
(hence $(C^2-C\cdot K_X)/2\ge0$ by Riemann-Roch). 
Assume $H\cdot C>0$.
Applying Riemann-Roch for surfaces
now gives $h^0(X,H+C)=h^0(X,H)+(C^2-C\cdot K_X)/2+H\cdot C>h^0(X,H)$.
Similarly, if $C$ is instead an exceptional curve with $C\cdot H>0$,
then $H+C$ is nef and effective, and $h^0(X,H+C)>h^0(X,H)$.
In particular, if $H$ is nef and effective with $|H|$
base curve free, and if $C$ is a prime divisor 
which is a base curve of $|H+C|$ such that $C$ is either exceptional
or $C^2>0$ (hence nef), then $C\cdot H=0$.

Now consider $Z=m_1P_1+\cdots+m_nP_n$ and let
$\alpha=\alpha(Z)$. Then $h_Z(t)=h^0(X,F)$, where
$F=tL-m_1E_1-\cdots-m_nE_n$. Let $H$ be the moving part
of $|F|$, and decompose the fixed part as $D+N$, where $D$ 
is the sum of the curves in the base locus of $|F|$
of nonnegative self-intersection and $N$
is the sum of the curves in the base locus
of negative self-intersection. 
By Conjecture \ref*{SHGHconjVanc}, each 
curve in $N$ is an exceptional curve, hence
disjoint from $H$ and $D$. In addition,
these exceptional curves are pairwise orthogonal 
(since if $C$ and $C'$ both appear in $N$
and have $C\cdot C'>0$, then $h^0(X, C+C')\ge C\cdot C'+1$
by Riemann-Roch, hence $C+C'$ cannot be part of the base locus of $|F|$).
Also, none of the exceptional curves $E$
appearing in $N$ is the blow up of a point
$P_i$ (since if it were we would have
$m_i=E\cdot F=E\cdot(H+D+N)=E\cdot N<0$).
Thus for some $c_{j_k}$ we have 
$N=\sum_kc_{j_k}C_{j_k}$, where $C_j$ is the proper transform
of the plane curve defined by $f_j$ in the statement of
Conjecture \ref*{simpleSHGH2}.
Applying Riemann-Roch gives
$$\begin{array}{lll}
h_Z(t) & = & h^0(X,F)=h^0(X,F-N)=((F-N)^2-K_X\cdot(F-N))/2+1\\ 
       & = & (F^2-K_XF)/2+1+\sum_k{c_{j_k}\choose 2}
={t+2\choose 2}-\sum_i{m_i+1\choose 2}+\sum_k{c_{j_k}\choose 2},
\end{array}$$
as claimed in Conjecture \ref*{simpleSHGH2}.

Now consider degree $t+1$, so $t+1>\alpha(Z)$.
Then $h_Z(t+1)=h^0(X, F+L)$, where 
$F+L=H+(L+D)+N$. Note that, as we saw above, 
$|L+D|$ is fixed component free. 
Thus the fixed part of $|F+L|$ consists
at most of exceptional curves coming from $N=\sum_kc_{j_k}C_{j_k}$. 
In fact, by the first paragraph of the proof,
$|F+D+L+\sum_k \hbox{min}(d_{j_k},c_{j_k})C_{j_k}|$ 
is base curve free and
$|F+D+L+N|=|F+D+L+\sum_k \hbox{min}(d_{j_k},c_{j_k})C_{j_k}|
+\sum_k \hbox{max}(0,c_{j_k}-d_{j_k})C_{j_k}$.
I.e., $\sum_k \hbox{max}(0,c_{j_k}-d_{j_k})C_{j_k}$
is the divisorial part of the base locus of $|F+L|$.
If we denote $\hbox{max}(0,c_{j_k}-d_{j_k})$ by $b_{j_k}$
and reindex, this becomes $\sum_l b_lC_l$, and applying Riemann-Roch as above
gives
$$h_Z(t+1)=h^0(X,F+L)={(t+1)+2\choose 2}-\sum_i{m_i+1\choose 2}+\sum_l{b_l\choose 2},$$
as claimed in Conjecture \ref*{simpleSHGH1}. (Note that $t+1$ here is the same as
$t$ in the statement of Conjecture \ref*{simpleSHGH1}, since there
we assumed $t>\alpha(Z)$,
but here, in order to handle Conjecture \ref*{simpleSHGH2}
simultaneously, we assumed only $t\ge\alpha(Z)$.) 
\prfend

\subsection{Mumford's Snake Lemma}\label{msl}

Mumford \cite{refMu1} applied the snake lemma
to questions related to $\mu_F$. We recall that now.
To do so we establish some notation that
we will use here and throughout the paper.
Let $F$, $C$ and $D$ be divisors on $X$ with $C$ effective; then 
we have the natural multiplication maps
$$\mu_{F,D}: H^0(X, F)\otimes H^0(X,D) \to H^0(X, F+D);$$
and
$$\mu_{C;F,D}: H^0(C, F\vert_C)\otimes H^0(X,D) 
\to H^0(C, (F+D)\vert_C).$$
In the particular case that $D=L$, which is almost always true in the
present paper, we write $\mu_{F}$ and $\mu_{C;F}$
instead of $\mu_{F,L}$ and $\mu_{C;F,L}$.

\begin{lem}\label{snakelem}
Let $p:X\to \pr2$ be a blow up of $\pr2$ at $n$ distinct points
with $L,E_1,\ldots,E_n$ as usual. Let $D$ be a divisor on $X$,
let $V=H^0(X, D)$, and let $F$ and $C$ be divisors on $X$ with $C$ effective
and with $h^1(X, F)=0=h^1(X, F+D)$. Then the following diagram
is commutative with exact rows:
{
$$\matrix {
0 & \to & H^0(X, F)\otimes V  & \to & H^0(X,F+C)\otimes V &
\to & H^0(C, (F+C)\vert_C)\otimes V & \to & 0 \cr
{} & {} & \downarrow \raise3pt\hbox to0in{$\scriptstyle\mu_{F,D}$\hss} & {}
& \downarrow\raise3pt\hbox to0in{$\scriptstyle\mu_{F+C,D}$\hss} & {}
& \downarrow\raise3pt\hbox to0in{$\scriptstyle\mu _{C;F+C,D}$\hss} &
{} & {} \cr
0 & \to & H^0(X,F+D) & \to & H^0(X,F+C+D) & \to &
H^0(C,(F+C+D)\vert_C)& \to & 0 \cr
} \eqno{(^\circ)}$$}
The snake lemma thus gives an exact sequence
$$\begin{array}{lllllll}
0 & \to & \hbox{ker}\,\mu_{F,D}  & \to & \hbox{ker}\,\mu_{F+C,D}  & \to &
\hbox{ker}\,\mu_{C;F+C,D} \\
   {}  & \to & \hbox{cok}\,\mu_{F,D}  & \to &  \hbox{cok}\,\mu_{F+C,D}  & \to
& \hbox{cok}\,\mu_{C;F+C,D}\to0
\end{array}$$
which we will refer to as ${\cal S}(F,C,D)$ (or ${\cal S}(F,C)$ if $D=L$).
\end{lem}

Another useful fact is the Castelnuovo-Mumford Lemma \cite{refMu2}, which
gives a criterion for $\mu_F$ to not only have maximal rank but to
be surjective. The version we state, Corollary \ref*{CMlem}, follows easily from
Lemma \ref*{snakelem}, using ${\cal S}(0,L,H)$.
(The hypothesis $h^1(X,H-L)=0$ is used to ensure that
$H^0(X, H)\to H^0(L,H|_L)$ is surjective; $H\cdot L\ge0$
then ensures that $\hbox{cok}\,\mu_{L;L,H}=0$.)

\begin{cor}\label{CMlem}
Let $p:X\to \pr2$ be obtained by blowing up $n$ distinct points of $\pr2$,
with $L$ the pullback of the class of a line.
If $H$ is a divisor on $X$ with $h^1(X,H-L)=0$ and $H\cdot L\ge0$, then
$\hbox{cok}\,\mu_{H}=0.$
\end{cor}

When $D=L$, explicit expressions for the kernels and cokernels in Lemma \ref*{snakelem}
can be given in terms of the cotangent bundle.
Recall the Euler sequence defining $\Omega =\Omega_{\pr2}$:
$$0 \to \Omega (1) \to {\cal O}_{\pr2} \otimes H^0(\pr2, {\cal O}_{\pr2}(1)) \to
{\cal O}_{\pr2}(1) \to 0.$$
Pulling the Euler sequence back to $X$, tensoring by ${\cal O}_X(F)$ and
identifying
$H^0(\pr2, {\cal O}_{\pr2}(1))$ with $H^0(X, L)$ gives an exact sequence
$$0 \to (p^*\Omega)(F+L) \to {\cal O}_{X}(F) \otimes H^0(X, L) \to {\cal
O}_{X}(F+L) \to 0. \eqno{(\dagger)}$$
If $h^1(X, F)=0$, then taking cohomology gives an exact sequence
$$0 \to H^0(X, (p^*\Omega)(F+L)) \to H^0(X, F)\otimes H^0(X, L) \to
H^0(X, F+L) \to H^1(X, (p^*\Omega)(F+L))\to0,$$
hence $\hbox{ker}\,\mu_{F}=H^0(X, (p^*\Omega)(F+L))$ and
$\hbox{cok}\,\mu_{F} =H^1(X, (p^*\Omega)(F+L))$.
Similarly, if $h^1(C, F|_C)=0$, by
restricting ($\dagger$) to $C$ and taking cohomology we see
$\hbox{ker}\,\mu_{C; F}=H^0(C, ((p^*\Omega)(F+L))|_C)$ and
$\hbox{cok}\,\mu_{C; F}=H^1(C, ((p^*\Omega)(F+L))|_C)$.
In case $C$ is a smooth rational curve, taking $t=F\cdot C$ we have:
$$\begin{array}{lll}
\hbox{ker}\,\mu_{C; F} & = & H^0(C, {\cal O}_C(t-a_C)\oplus{\cal O}_C(t-b_C))\\
\hbox{cok}\,\mu_{C; F} & = & H^1(C, {\cal O}_C(t-a_C)\oplus{\cal O}_C(t-b_C))\\
\end{array}\eqno{(\dagger\dagger)}$$
with $a_C$ and $b_C$ as defined in section \ref*{bnc}.

\renewcommand{\thethm}{\thesection.\arabic{thm}}
\setcounter{thm}{0}

\section{Main Results}\label{mainresults}

Given a divisor $F$ on a blow up $X$ of $\pr2$ at $n$ general points,
the naive conjecture that $h^0(X, F)$ always equals
$\hbox{max}(0, (F^2-K_X\cdot F)/2+1)$ is false.
One way to salvage it, is to impose a niceness requirement on $F$,
such as to require that $F$ be nef, or even
that $F\cdot C\ge 0$ for all exceptional $C$, which is weaker. In fact, 
an equivalent version of the SHGH Conjecture is given by Conjecture \ref*{SHGHconj},
which just conjectures 
that $h^0(X, F)=\hbox{max}(0, (F^2-K_X\cdot F)/2+1)$
whenever $F\cdot C\ge0$ for all exceptional $C$.

Likewise, the naive conjecture that $\mu_F$ should always have maximal rank
is also false, as we saw above.
Again we can try to salvage the naive conjecture by 
imposing a niceness requirement on $F$, but the necessary requirement is more subtle.
Assuming we can compute $h^0(X, F)$ for an arbitrary divisor $F$, 
Lemma \ref*{baselocuslem} reduces the problem of computing the rank
of $\mu_F$ in general to the case that $F$ is effective and fixed component free.
But, as we mentioned above,
even if $F$ is effective and fixed component free, or even very ample,
$\mu_F$ can fail to have maximal rank. 

Instead, we will consider all $F\in L+\hbox{EFF}(X)$,
where $\hbox{EFF}(X)$ is the subsemigroup of 
the divisor class group of $X$ of classes of all effective divisors.
In this section, refining ideas of Fitchett, we will
give an upper bound on the dimension of the cokernel of
$\mu_F$ for certain $F\in L+\hbox{EFF}(X)$
(for all of them if the SHGH Conjecture is true).
We conjecture that this upper bound is in fact an equality.

Assuming the SHGH Conjecture, Fitchett reduced the problem
of handling $\mu_F$ for an arbitrary $F\in L+\hbox{EFF}(X)$
to the case $F=L+mE$ where $E$ is an exceptional curve
and $0\le m\le L\cdot E$.
(See \cite{refF1}, \cite{refF2}, which give explicit bounds
on $\hbox{dim cok}\,\mu_F$ for $F\in L+\hbox{EFF}$ in the case
of $n\le8$ general points, using a construction
originally described in Fitchett's thesis.)

We now recall Fitchett's idea, assuming that $X$ is obtained
by blowing up $n$ general points. Note that the SHGH Conjecture
would make the assumption $h^1(X, H)=0$ automatic.

\begin{prop}\label{fitchettidea}
Let $F\in L+\hbox{EFF}(X)$, so we have 
the decomposition $F-L=H+N$ given by Lemma \ref*{Wlem}(d,e), 
where $H$ is effective with $H\cdot E\ge0$
for all exceptional $E$ and where either $N=0$ or
$N=c_1C_1+\cdots+c_rC_r$ for some mutually 
disjoint exceptional curves $C_i$
and integers $c_i>0$. Assume that $h^1(X, H)=0$. 
Then, $cok\,\mu_F\cong \oplus _{i=1,\dots,r}
cok\,\mu_{L+c_iC_i}\cong cok\,\mu_{L+N}$. If moreover $F\cdot C_i\ge0$, then
$0< c_i\le L\cdot C_i$.
\end{prop}

\noindent{\bf Proof.} By Lemma \ref*{Wlem}(d,e), $F=L+H+N$ where either $N=0$ or
$N=c_1C_1+\cdots+c_rC_r$ for some mutually disjoint exceptional curves $C_i$
and integers $c_i>0$, and where $H$ is effective
and orthogonal to $N$.
Then $h^1(X, H)=0$ implies $cok\,\mu_{L+H}=0$ by Corollary \ref*{CMlem}.
Taking cohomology of
$$0\to {\cal O}_X(H)\to {\cal O}_X(L+H)\to {\cal O}_L(H+L)\to 0$$
and using  $h^1(X, H)=0$ implies that $h^1(X, L+H)=0$.
Similarly, we also have $h^1(X,2L+H)=0$.
So sequence ${\cal S}(F-N,N)$ of Lemma \ref*{snakelem} holds (since $F-N=L+H$) and tells
us that
the cokernels for  $\mu_F$ and $\mu_{N;F}$
are isomorphic. But ${\cal O}_N(F)$ is isomorphic to
${\cal O}_N(L+N)$, since $H\cdot N=0$, and the
$C_i$ are disjoint, so ${\cal O}_N(L+N)\cong \oplus {\cal
O}_{c_iC_i}(L+c_iC_i)$; we finally get that the cokernel of
$\mu_{N;F}$ is isomorphic to the direct sum of the
cokernels of $\mu_{c_iC_i;L+c_iC_i}$. Moreover, $h^1(X,L)=h^1(X,2L)=0$,
so sequence ${\cal S}(L,c_iC_i)$ of Lemma \ref*{snakelem} gives
$cok\,\mu_{c_iC_i;L+c_iC_i} \cong cok\,\mu_{L+c_iC_i}$ since
$\hbox{cok}\,\mu_L=0$.
Hence, as Fitchett observed, the cokernel for  $\mu_F$
is isomorphic to the direct sum of the cokernels
for $\mu_{L+c_iC_i}$. Running the same argument with
$L+N$ in place of $F$ now gives $\oplus _{i=1,\dots,r}
cok\,\mu_{L+c_iC_i}\cong \mu_{L+N}$.\prfend

This and Lemma \ref*{baselocuslem} motivate the following problem:  

\begin{prob}\label{multexcdiv}  Determine the rank of $\mu_F$ for each $F=L+iE$, where
 $E$ is smooth and rational with $E^2=-1$ and $0\le i\le E\cdot L $.
\end{prob}

Concerning this problem, we prove Theorem \ref*{fitchettbounds}, which  
gives explicit upper bounds for the dimension of $\hbox{cok}\,\mu_{L+iE}$.
Similar but less precise results 
were given in \cite{refF1}, \cite{refF2}. Theorem \ref*{fitchettbounds}(c)
seems to be entirely new, however.

One can also give a lower bound for $\hbox{dim cok}\,\mu_F$. This lower bound,
as is the case for the upper bound given in Theorem \ref*{fitchettbounds},
is just what one can conclude from the sequence of short
exact sequences $(\ddagger)$. An explicit formula for this lower bound
turns out to be more complicated and less useful than that for
the upper bound, so we do not include it here.

\begin{thm}\label{fitchettbounds} Let $X$ be the blow up of
$\pr2$ at $n$ distinct points $P_1,\ldots,P_n$, and take
$L,E_1,\ldots,E_n$ as usual. Let $F=L+iE$, where
$d=E\cdot L$, $m$ is the maximum of $E\cdot E_j$ over
$1\le j\le n$, and where $E$
is smooth and rational with $E^2=-1$ and $0\le i\le d$. 
Then we have $$\hbox{dim cok}\,\mu_F\le {i-b_E\choose 2} + {i-a_E\choose 2},$$
with equality in the following cases:
\begin{description}
\item[(a)] $i\le a_E+2$;
\item[(b)] $b_E-a_E\le 2$; or
\item[(c)] $a_E=d-m$.
\end{description}
\end{thm}

\noindent{\bf Proof.} Let $F_j=L+jE$ for $0\le j<i$; it is easy to check that
$H^1(F_j)=0=H^1(F_j+L)$, so we can consider  ${\cal S}(F_j,E)$. 
Since $(F_j+E)\cdot E=d-j\ge0$, by $(\dagger\dagger)$ 
and Lemma \ref*{splitlem} we have 
$$\hbox{ker}\,\mu_{E;F_j} = H^0(E,{\cal O}_E(b_E-j))\oplus
H^0(E,{\cal O}_C(a_E-j)),$$
$\hbox{cok}\,\mu_{E;F_j} = H^1(E,{\cal O}_E(b_E-j))\oplus
H^1(E,{\cal O}_E(a_E-j))$, and
$\hbox{dim ker}\,\mu_{E;F_j}=(b_E-j+1)_+ + (a_E-j+1)_+$,
while $\hbox{dim cok}\,\mu_{E;F_j}=(j-b_E-1)_+ + (j-a_E-1)_+$.

Writing $a$ and $b$ for $a_E$ and $b_E$, we have the following
exact sequences:

$$\begin{array}{lllll}
\hbox{cok}\,\mu_{F_{i-1}}  & \to & \hbox{cok}\,\mu_{F_i}  &
\to &  H^1({\cal O}_E(b-i)\oplus {\cal O}_E(a-i)) \to  0\\
\hbox{cok}\,\mu_{F_{i-2}}  & \to & \hbox{cok}\,\mu_{F_{i-1}}  &
\to &  H^1({\cal O}_E(b-i+1)\oplus {\cal O}_E(a-i+1)) \to  0\\
\cdots & & & &   \\
\hbox{cok}\,\mu_{L+E}  & \to & \hbox{cok}\,\mu_{L+2E}  &
\to &  H^1({\cal O}_E(b-2)\oplus {\cal O}_E(a-2)) \to  0\\
\hbox{cok}\,\mu_{L}  & \to & \hbox{cok}\,\mu_{L+E}  &
\to &  H^1({\cal O}_E(b-1)\oplus H^1({\cal O}_E(a-1))  \to  0
\end{array}\eqno(\ddagger)$$

Note that $\hbox{cok}\,\mu_{L} =0$; this is just the fact that
$R_1\otimes R_1$ maps by multiplication surjectively to $R_2$,
where $R$ is the ring $K[\pr2]$.
Since $\hbox{cok}\,\mu_{L} =0$,
$\hbox{dim cok}\,\mu_{F_i}$ is at most
the sum of the dimensions of the column of $H^1$'s; i.e.,
we have $\hbox{dim cok}\,\mu_{F_i}\le \sum_{j\ge0} ((i-b-1-j)_+ + (i-a-1-j)_+)
= {i-b\choose 2} + {i-a\choose 2}$,
and equality holds if and only if
the displayed sequences are all exact
on the left. Moreover,
$h^1({\cal O}_E(b-j)\oplus {\cal O}_E(a-j))=0$
for all $j\le a+1$, so $\hbox{cok}\,\mu_{F_{j}}=0$
for all $0\le j\le a+1$. Thus each sequence is exact for which
$F_j$ in the middle column has index $j\le a+2$.
This implies claim (a).
Moreover, any of the sequences for which
$\hbox{dim ker}\,\mu_{E;F_j}=0$
will also be exact on cokernels, and
$\hbox{dim ker}\,\mu_{E;F_j}=0$ for all
$j\ge b+1$. It follows that equality holds if $j>a+2$
implies $j>b$; i.e., if $a+2\ge b$. This shows (b).

Finally, consider (c); thus $b=m$. 
It is enough to show that the maps
$\hbox{ker}\,\mu_{F_{j}}  \to \hbox{ker}\,\mu_{E;F_{j}}$
are onto for $a+3\le j\le b$; we already observed above that
exactness holds on cokernels (and hence for kernels)
for other values of $j$. We may assume, after reindexing
if need be, that $E_1\cdot E=m$.

From the exact sequence
$0\to {\cal O}_E(m-d)
\to {\cal O}_E\otimes H^0(X, L-E_1) \to {\cal O}_E(d-m) \to 0$
we see that $\hbox{ker}\,\mu_{E;F_{j},L-E_1}\cong 
H^0(E,F_j\cdot E+m-d)=H^0(E,m-j)$. 
The inclusion $H^0(X, L-E_1)\subset H^0(X, L)$
induces an inclusion $\hbox{ker}\,\mu_{E;F_{j},L-E_1}\to 
\hbox{ker}\,\mu_{E;F_{j},L}=\hbox{ker}\,\mu_{E;F_{j}}$.
The cokernel is isomorphic to $H^0(E, {\cal O}_E(a-j))$
(see exact sequence (4) in the proof of Theorem 3.1 of \cite{refF1}).
Unwinding definitions, we see that the induced map
$H^0(E,m-j)\cong \hbox{ker}\,\mu_{E;F_{j},L-E_1}
\hookrightarrow \hbox{ker}\,\mu_{E;F_{j}}$
sends an element 
$\sigma\in H^0(E,m-j)$ to 
$x|_E\sigma\otimes y-y|_E\sigma\otimes x\in \hbox{ker}\,\mu_{E;F_{j}}$,
if we choose homogeneous coordinates
$x,y$ and $z$ on $\pr2$ such that $P_1$ is the point where $x=0=y$.

Since $a+3\le j\le b=m$, the induced inclusion
$H^0(E,m-j)\hookrightarrow \hbox{ker}\,\mu_{E;F_{j}}$
is an isomorphism.
Moreover, taking cohomology of
$0\to {\cal O}_X(E_1+(j-1)E)
\to {\cal O}_X(E_1+jE) \to {\cal O}_E(m-j) \to 0$
gives the map $H^0(X, E_1+jE)\to H^0(E,{\cal O}_E(m-j))$,
which is surjective for $1\le j\le m+1$
since $h^1(E, {\cal O}_E(m-j))$
is 0 in this range and, by induction on $j$
starting with $j=1$, so is
$h^1(X, {\cal O}_X(E_1+(j-1)E))$.
Composing $H^0(X, E_1+jE)\to H^0(E,{\cal O}_E(m-j))$
with the induced isomorphism
$H^0(E,m-j)\to \hbox{ker}\,\mu_{E;F_{j}}$
gives for each $f\in H^0(X, E_1+jE)$
the map $f\mapsto f|_E\mapsto x|_Ef|_E\otimes y-y|_Ef|_E\otimes x
\in \hbox{ker}\,\mu_{E;F_{j}}$. Thus every element of
$\hbox{ker}\,\mu_{E;F_{j}}$ is of the form
$f|_E\mapsto x|_Ef|_E\otimes y-y|_Ef|_E\otimes x$ 
where $f\in H^0(X, E_1+jE)$. But
$f|_E\mapsto x|_Ef|_E\otimes y-y|_Ef|_E\otimes x$ is the image
of $xf\otimes y-yf\otimes x\in \hbox{ker}\,\mu_{F_{j}}$
under the map $\hbox{ker}\,\mu_{F_{j}}  \to \hbox{ker}\,\mu_{E;F_{j}}$,
so the map is surjective.
\prfend

In fact, we do not know any times that the equality in the theorem
does not hold. This suggests the following conjecture:

\begin{conj}\label{boundseqconj} Let $F=L+iE$ be as in Theorem
\ref*{fitchettbounds} for general points $P_i$.
Then $\hbox{dim}\,\hbox{cok}\,\mu_F= {i-b_E\choose 2} + {i-a_E\choose 2}$.
\end{conj}

\begin{remark}\label{betticonjrem}\rm
(a) Conjecture \ref*{boundseqconj} is equivalent to the first column of maps
in $(\ddagger)$ all being injective. Hence if Conjecture \ref*{boundseqconj}
holds for some $i$, then it holds for all $0\le j\le i$. Moreover, 
Conjecture \ref*{boundseqconj}
holds for $i=L\cdot E$ if and only if it holds for all
$a_E+3\le i\le b_E$, since the proof of Theorem \ref*{fitchettbounds}
shows that the first column of maps in
$(\ddagger)$ are injective for 
$i\le a_E+2$ and for $i>b_E$.

(b) In the notations of Proposition \ref*{fitchettidea}, assuming the SHGH Conjecture and
Conjecture \ref*{boundseqconj}, and assuming we can determine splitting types, we thus can determine the
dimension of the cokernel of $\mu_F$ for any $F$ as long as
$|F-L|$ is not empty. The splitting type of an exceptional curve
can be computed fairly efficiently, at least provisionally
(that is, by Macaulay 2 \cite{refGS}, say, in positive
characteristic, using randomly
chosen points; a Macaulay 2 script that does this is included
in Section A2.3 of the posted version
of this paper, \cite{refGHIpv}). We discuss this,
and we give additional, computational, support for Conjecture \ref*{boundseqconj},
in Section \ref*{algrthms} of the appendix.

(c) Translating in terms of fat points, this says that we can produce
conjectural dimensions for the cokernels of $\mu_t$ for $I(Z)$
in every degree $t$ but $t=\alpha(Z)$. We can even sometimes
determine the dimension for the cokernels of $\mu_{\alpha}$, for example 
by applying Lemma \ref*{qstarlstarlem} (see Example \ref*{boundsexample}),
or if $h_Z(\alpha)=1$, or if $F_{\alpha}(Z)$ decomposes as $F_{\alpha}(Z)=L+H+N$ 
where $h^1(X, H)=0$ and $H\cdot L\ge0$, even if $H$ is not effective (see Example \ref*{rarecases}).

(d) Conjecture \ref*{Betticonj} is equivalent to Conjecture \ref*{boundseqconj},
assuming the SHGH Conjecture. The first sum in the bound
in Theorem \ref*{Bettibound} is exactly the difference term
in Lemma \ref*{baselocuslem}, which accounts for the contribution to the
cokernel owing to loss in fixed components 
in going from degree $t$ to degree $t+1$. (This term does not occur
in Conjecture \ref*{boundseqconj} since the $F$ there is base curve free.)
The second sum in the bound in Theorem \ref*{Bettibound} sums up exactly 
what each of the disjoint exceptional curves 
in the base locus of $F_t(Z)-L$ should contribute 
to the cokernel, according to 
Proposition \ref*{fitchettidea} and Conjecture \ref*{boundseqconj}.
\end{remark}

\noindent{\it Proof of Theorem \ref*{Bettibound}}:
Let $t>\alpha(Z)$ and let $F=F_t(Z)$.
Thus $F-L$ is effective. 
By the SHGH Conjecture (\ref*{SHGHconjVanc}), as 
in the proof of Theorem \ref*{conjequiv},
we have $F-L=H+N$, where $H$ is nef and effective
with $h^1(X, H)=0$, and where
$N=c_1C_1+\cdots+c_rC_r$ is a sum of pairwise orthogonal 
exceptional curves orthogonal to $H$
with the curves $C_j$ being the proper 
transforms of those curves in the base locus of
$I(Z)_{t-1}$ which are negative for the points $P_i$.
Let $N'$ be that part of $N$
which remains in the base locus for $|F|$ and let $N''$ be 
what remains in the base locus of $|F+L|$;
thus $N'=c_1'C_1+\cdots+c_r'C_r$ and $N''=c_1''C_1+\cdots+c_r''C_r$.
Note that $c_j-c'_j\leq \hbox{deg}(C_j)$, because
$-c_j=F_{t-1}(Z)\cdot C_j$ and $-c_j\le-c'_j=
\hbox{min}(F_{t}(Z)\cdot C_j,0)
\le L\cdot C_j + F_{t-1}(Z)\cdot C_j = \hbox{deg}(C_j)-c_j$.
By Lemma \ref*{baselocuslem},
$\hbox{dim cok }\mu_t=\hbox{dim cok }\mu_F=
\hbox{dim cok }\mu_{H+L+N-N'}+(h^0(X, F+L)-h^0(X, H+2L+N-N'))$.
By Proposition \ref*{fitchettidea} and Theorem \ref*{fitchettbounds},
$\hbox{dim cok }\mu_{H+L+N-N'}=
\sum_j\hbox{dim cok }\mu_{L+(c_j-c_j')N_j}\le 
\sum_j\Bigg({c_j-c_j'-a_{C_j}\choose 2}
+{c_j-c_j'-b_{C_j}\choose 2}\Bigg)$. 

Since $N''$ is in the base locus of $|F+L|$, we have
$h^0(X, F+L)=h^0(X, F+L-N'')$. But $F+L-N''=H+2L+N-N''=H+2L+N-N'+(N'-N'')$,
and both $H+2L+N-N''$ and $H+2L+N-N'$ are nef and effective, so
by Conjecture \ref*{SHGHconjVanc} we have
$h^1(X, H+2L+N-N'+(N'-N''))=0=h^1(X,H+2L+N-N')$.
Plugging into Riemann-Roch and simplifying gives $h^0(X, F+L)-h^0(X, H+2L+N-N')=
h^0(X, H+2L+N-N'+(N'-N''))-h^0(X, H+2L+N-N')=
(N'-N'')^2/2-K_X\cdot(N'-N'')/2 + (2L+H+N-N')\cdot(N'-N'')$.
Keeping in mind that $(L+H+N-N')\cdot(N'-N'')=0$, this gives
$\sum_jd_j(c_j'-c_j'')-\sum_j{c_j'-c_j''\choose 2}$. 
Putting everything together gives 
$\hbox{dim cok }\mu_t=\le 
\sum_j\Bigg({c_j-c_j'-a_{C_j}\choose 2}
+{c_j-c_j'-b_{C_j}\choose 2}\Bigg)+\sum_jd_j(c_j'-c_j'')-\sum_j{c_j'-c_j''\choose 2}$.
\prfend

\newcounter{myapp} 
\newcounter{mythm} 
\setcounter{myapp}{0}
\renewcommand{\thesection}{Appendix}
\renewcommand{\thesubsection}{A\arabic{section}.\arabic{subsection}}
\setcounter{mythm}{0}
\setcounter{section}{0}
\renewcommand{\thethm}{A\arabic{section}.\arabic{subsection}.\arabic{mythm}}

\vfil\eject

\noindent{\LARGE \bf Appendix}

\setcounter{section}{0}
\renewcommand{\thesection}{A\arabic{section}}
\renewcommand{\thesubsection}{A\arabic{section}.\arabic{subsection}}

\addtocounter{myapp}{1}
\section{Making the SHGH Conjecture Explicit}\label{explicitSHGH}

In this appendix we give another version of the
SHGH Conjecture and show how to derive explicit predictions
for values of Hilbert functions using it.

\setcounter{mythm}{0}
\subsection{The Weyl Group}\label{TWG}

We now recall the Weyl group $W=W_n$, which acts on $\hbox{Cl}(X)$
but which depends only on
the number $n$ of points $P_i$ blown up.  If $0\leq n\leq1$, then $W=\{id\}$ is trivial.
If $n=2$, then $W=\{id,s_1\}$, where for any divisor class $F$,
$s_1(F)=F+(r_1\cdot F)r_1$, where $r_1=E_1-E_2$.
For $n>2$, let $r_0=L-E_1-E_2-E_3$ and for $1\leq i<n$, let
$r_i=E_i-E_{i+1}$. Then $W$ is generated by the operators
$s_i$, $0\leq i<n$, where $s_i(F)=F+(r_i\cdot F)r_i$.
It is now easy to check that $W$ preserves the intersection form
(i.e., $wF\cdot wG
=F\cdot G$ for all $w\in W$ and all $F,G\in \hbox{Cl}(X)$), and that
$wK_X=K_X$ for all $w\in W$. The subgroup generated by 
$s_1,\ldots, s_{n-1}$ is just the permutation group on
$E_1,\ldots,E_n$. The action of the element $s_0$
corresponds to that of the quadratic Cremona transformation
centered at $P_1,P_2,P_3$.

If $n=0$, $X=\pr2$ has no exceptional curves.
If $n=1$, then $E_1$ is the only exceptional curve, and if
$n=2$, then $E_1$, $E_2$ and $L-E_1-E_2$ are the only exceptional curves.
For $n\geq3$,  Nagata \cite{refN2} has shown that if a class
$E$ is the class of an exceptional curve, then $E\in WE_n$; i.e., 
the classes of exceptional curves lie in a single $W$-orbit.
Moreover, if $E\in WE_n$ and the points $P_i$ are general, Nagata showed 
$E$ is the class of an exceptional curve. (When $E\in WE_n$
but the points are not general, then although $E$ is effective,
it can fail to be reduced and irreducible, and thus need not be
the class of an exceptional curve. For example, $2L-E_1-\cdots-E_5\in WE_5$,
but if $P_1,P_2,P_3$ are collinear, then 
the proper transform of the line through $P_1,P_2,P_3$
is a fixed component of $|2L-E_1-\cdots-E_5|$.)

If $n=0$, let ${\cal E}_n$ be the submonoid of 
$\hbox{Cl}(X)$ generated by $L$. If $n=1$, let ${\cal E}_n$ 
be the submonoid generated by $L-E_1$ and $E_1$. If 
$n=2$, let ${\cal E}_n$ be generated by $L-E_1-E_2$, $E_1$ and $E_2$,
while if $n\ge3$, let ${\cal E}_n$ be generated by $L-E_1-E_2$ and $E_1,\ldots,E_n$
(or equivalently by the orbit $WE_n$ of $E_n$ under $W$).
Thus if
$n\geq 3$ every element $D\in {\cal E}_n$ is of the form $D=\sum c_iC_i$,
where $c_i$ is a nonnegative integer and $C_i$ is an exceptional curve.
Define ${\cal E}^*_n$ to be the dual cone;
thus $F\in{\cal E}^*_n$ means that $F\cdot D\geq0$ for every $D\in{\cal E}_n$
(and thus that $F\cdot E\geq0$ for every exceptional curve
$E$).

Let $\hbox{EFF}=\hbox{EFF}(X)\subset \hbox{Cl}(X)$ 
denote the submonoid of classes of effective divisors, 
let $\hbox{NEF}=\hbox{NEF}(X)\subset \hbox{Cl}(X)$ denote the
submonoid (indeed the cone, since a class is nef if a positive multiple
is) of classes of nef divisors, let $\Psi_n$ be the submonoid 
generated by the union of ${\cal E}_n$ and the element $-K_X$
and let $\Delta_n$ be the submonoid of $\hbox{Cl}(X)$
generated by $H_0=L$, $H_1=L-E_1$, $H_2=2L-E_1-E_2$, and
$H_i=-K_X+E_{i+1}+\cdots+E_n=3L-E_1-\cdots-E_i$, for $3\leq i\leq n$.
Notice that  $F=tL-m_1E_1-\cdots-m_nE_n\in\Delta_n$
if and only if $t\ge m_1+m_2+m_3$ and
$m_1\geq m_2\geq \cdots \geq m_n\geq0$, and 
that $H_n$ is $-K_X$ and $H_i\cdot H_j\geq 0$ unless $i,j\geq
10$. 

Since the next result is a statement for all classes on $X$, we need
to state it in terms of a blowing up of generic points. However,
when we are interested in a specific class $F$, it is enough
to consider a blow up of general points (but the conditions
of generality will depend on $F$). The following result is known
but hard to cite.

\addtocounter{mythm}{1}
\begin{lem}\label{Wlem}
Let $X$ be the blow up of $\pr2$ at $n$ generic points $P_i$.
\begin{description}
\item[(a)] If $A\in\Delta_n$ and $w\in W_n(X)$, then $wA=
A+a_0r_0+\cdots+a_{n-1}r_{n-1}$ for some nonnegative integers $a_i$.
\item[(b)] ${\rm NEF}(X)\subset {\cal E}^*_n= W_n \Delta_n\subset \Psi_n$
\item[(c)] $h^j(X,F)=h^j(X,wF)$ for all $j$, all $w\in W_n$ and all $F\in \hbox{Cl}(X)$
\item[(d)] ${\rm EFF}(X)\subset \Psi_n$
\item[(e)] If $F\in\Psi_n$, then there is a unique decomposition
$F=H+N$ where $H\in {\cal E}^*_n$, $N\in{\cal E}_n$, $H\cdot N=0$
and, if $N\ne0$, then $N=c_1C_1+\cdots+c_rC_r$
where for each $i$, $c_i$ is a positive integer and $C_i$ is the class
of an exceptional curve with $C_i\cdot C_j=0$ for all $i\ne j$.
Moreover, $H$ is effective if $F$ is.
\end{description}
\end{lem}

\noindent {\em Proof}. (a) See, for example, Lemma 1.2 (1) of \cite{refMA}.
(b) First note that ${\rm NEF}\subset {\cal E}^*_n$, since
${\cal E}_n\subset{\rm EFF}$. 
Now we verify that ${\cal E}^*_n= W_n\Delta_n$. We leave the 
cases $0\le n<3$ to the reader; assume $n\ge 3$. It is known that
$\hbox{EFF}$ and $\hbox{NEF}$ are $W_n$ invariant,
and that the set of exceptional divisors is just the orbit
$W_nE_n$; cf. \cite{refN2}.
Since $W_n$ preserves the set of exceptional curves,
if we show $\Delta_n\subset{\cal E}^*_n$, then
$W_n\Delta_n\subset{\cal E}^*_n$. But $H_i$ is nef
(hence in ${\cal E}^*_n$) for $i\leq9$.
For $i>9$, $H_i=-K_X+E_{i+1}+\cdots+E_n$,
hence for any exceptional curve $E$ we have
$H_i\cdot E\geq0$, since $E$ meets $-K_X$ once
and $E\cdot E_j\geq-1$ with equality if and only if
$E=E_j$. Thus $\Delta_n\subset{\cal E}^*_n$.

Conversely, say $F\in {\cal E}^*_n$.
Note that every element $D\in{\cal E}^*_n$
satisfies $D\cdot L\geq0$, since
$L=(L-E_1-E_2)+E_1+E_2$. Since $W_n$ preserves
${\cal E}^*_n$, there must be some $w\in W_n$ such that
$L\cdot wF$ is as small as possible.
Thus, $wF\cdot r_0\geq0$, otherwise we would
have $s_0wF\cdot L<wF\cdot L$.
We can also assume
$m_1\geq m_2\geq \cdots \geq m_n$,  where
$m_i=wF\cdot E_i$, since each operator $s_i$, $i>0$,
merely transposes $E_i$ and $E_{i+1}$, so we can in $W_n$ permute
the $E_i$ without affecting $L\cdot wF$.
Thus $wF\cdot r_i\geq0$ for all $i\geq1$.
Finally, $wF\cdot (L-E_1-E_2)\geq0$ since $wF\in {\cal E}^*_n$, and
$wF\cdot (L-E_1)\geq0$ since $L-E_1=(L-E_1-E_2)+E_2$
is a sum of exceptional curves.
By Lemma 1.4 of \cite{refTrans}, we thus have $wF\in \Delta_n$.

Finally, it is clear that $\Delta_n\subset \Psi_n$
but $\Psi_n$ is $W_n$-invariant, so $W_n\Delta_n\subset\Psi_n$.

(c) This is, in somewhat different language, due to Nagata \cite{refN2}.
The basic idea is this. Say $F=tL-m_1E_1-\cdots-m_nE_n$.
Then $wF=twL-m_1wE_1-\cdots-m_nwE_n$, where $L'=wL,E'_1=wE_1,\ldots,E'_n=wE_n$
is a basis of $\hbox{Cl}(X)$. This basis is, however, an {\it exceptional configuration};
i.e., there is a birational morphism $p':X\to\pr2$
such that $E'_i=p'^{-1}(P'_i)$ for some points $P'_i\in\pr2$
and such that $L'$ is the pullback via $p'$ of 
the class of a line in $\pr2$ (see Theorem 0.1 of \cite{refDuke}), 
but the points $P'_i$ are themselves generic. 
Since $P_i$ and $P'_i$ both give generic sets of points, all that matters are
the coefficients $t$ and $m_i$, so we have
$h^j(X,tL-m_1E_1-\cdots-m_nE_n)=h^j(X,tL'-m_1E'_1-\cdots-m_nE'_n)$;
i.e., $h^j(X, F)=h^j(X,wF)$.

(d) Let $F\in\hbox{EFF}(X)$. Then there are at most finitely many exceptional
curves $E$ such that $F\cdot E<0$. Let this finite set of distinct exceptional
curves be $C_1,\cdots,C_r$, 
let $c_i=-F\cdot C_i$, let $N=c_1C_1+\cdots+c_rC_r$ and let $H=F-N$.
Note that $C_i\cdot C_j=0$ for all $i\ne j$ (since $(C_i+C_j)\cdot F<0$,
but $C_i\cdot C_j>0$ implies $(C_i+C_j)$ meets both $C_i$ and $C_j$
nonnegatively and hence is nef) and that $H\cdot C_i=0$. 
Since $F$ is effective, $N$ is contained in the scheme 
theoretic base locus of $|F|$, hence $H$ is effective.
But if $H\cdot E<0$, then $E$ is not $C_i$ for any $i$,
hence $E\cdot F\ge0$ so we get $E\cdot C_i>0$ (implying $E+C_i$ is nef)
for some $i$ even though $H\cdot(E+C_i)<0$.
It follows that $H\in {\cal E}^*_n$.
The result follows since $N\in {\cal E}_n\subset \Psi_n$ and 
$H\in{\cal E}^*_n\subset\Psi_n$. 

(e) We leave the cases $0\le n<3$ to the reader, so assume $n\ge3$.
If $F\in\Psi_n$, then $wF\cdot L\ge 0$ and $wF\cdot (L-E_1)\ge 0$ 
for all $w\in W_n$, since $wF\in \Psi_n$ but
$L$ and $L-E_1$ are nef and meet $-K_X$ nonnegatively.
Choose $w$ such that $wF\cdot L$ is as small as possible, and write
$wF=tL-m_1E_1-\cdots-m_nE_n$. Since $W_n$ includes the group of 
permutations of $E_1,\ldots,E_n$, we may assume  that
$m_1\ge \cdots\ge m_n$. Since $t$ is as small as possible,
we know that $t\ge m_1+m_2+m_3$ (otherwise $s_0wF\cdot L<wF\cdot L$).
If $m_3\ge0$ or $m_2\le0$, let $H'=tL-\sum_{m_i>0}m_iE_i$ and $N'=-\sum_{m_i<0}m_iE_i$. 
Using the definition of $\Delta_n$ and the facts that
$t\ge m_1+m_2+m_3$ and $t\ge m_1$, 
it is now not hard to check that $H'\in\Delta_n$, and clearly $N'\in{\cal E}_n$.
If $m_3<0$ and $m_2>0$, there are two cases. If $c=(tL-m_1E_1-m_2E_2)\cdot(L-E_1-E_2)<0$,
then let $H'=(t+c)L-(m_1+c)E_1-(m_2+c)E_2$, and let
$N'=(-c)(L-E_1-E_2)-m_3E_3-\cdots-m_nE_n$. If $(tL-m_1E_1-m_2E_2)\cdot(L-E_1-E_2)\ge0$,
let $H'=tL-m_1E_1-m_2E_2$, and let
$N'=-m_3E_3-\cdots-m_nE_n$. Either way $H'\in\Delta_n$ and $N'\in{\cal E}_n$.
In all cases we also have $H'\cdot N'=0$ and that the components of $N'$ are disjoint
and orthogonal to $H'$.
We now take $H=w^{-1}H'$ and $N=w^{-1}N'$, where the classes $C_i$
are the $w^{-1}$ translates of the components of $N'$.
Uniqueness follows from the fact that $N=\sum_{E\hbox{ exceptional}}-(F\cdot E)E$.
Note that in case $F$ is effective we found a decomposition $F=H+N$
in (d), with $H$ effective. Uniqueness now shows that $H$ is necessarily effective. 
\prfend

\addtocounter{mythm}{1}
\begin{remark}\label{excsuffice}\rm
As an application of Lemma \ref*{Wlem}, we will classify all smooth rational curves $C$
on a blow up $p:X\to \pr2$ of general points, either by assuming the SHGH Conjecture,
or by working over the complex numbers using Proposition 2.4 of \cite{refdF},
assuming that the points blown up are very general.
In either case, we have $C^2\geq -1$.
If $C^2=-1$, then $C$ is an exceptional curve.
If $C^2>-1$, then $C$ is nef, hence $wC\in\Delta$ for some $w\in W$ by Lemma \ref*{Wlem}(b),
so it suffices if we find all smooth rational $C\in\Delta$.
Since $C\in\Delta$, we have $C=\sum_i a_iH_i$ for some
nonnegative integers $a_i$.
Note that $C\cdot H_j\geq0$ for $0\leq j\leq2$.
By adjunction and $C^2>-1$ we have $-C\cdot K_X=C^2+2\ge 2$.
Since $C$ is nef, we have $C\cdot (E_{j+1}+\cdots+E_n)\ge0$.
Thus $C\cdot H_j=C\cdot (-K_X+E_{j+1}+\cdots+E_n)
\geq 2$ for all $j>2$.
If $a_j>0$ for some $j>2$, let $C'=C-H_j$. Since
$C'$ is still a nonnegative
integer combination of the $H_i$, we have $C\cdot C'\geq0$. Now
$C^2=C\cdot C' + C\cdot H_j\geq C\cdot H_j\geq C\cdot H_n=-C\cdot K_X=C^2+2$.
I.e., we must have $a_j=0$ for all $j>2$.

Thus $C=aH_0+bH_1+cH_2$. It is now an easy exercise using adjunction
to show that
the only solutions to $C^2+C\cdot K_X=-2$ are
$H_0$, $2H_0$, $H_1$, $H_2$, $H_0+bH_1$ and $H_2+bH_1$.
Thus $W$ orbits of these and $E_1$ are the only possible
smooth rational curves in $X$.

Each such $C$ can be turned into an exceptional curve $E$ by subtracting
off additional $E_i$; for example,
if $C=H_2+2H_1=4L-3E_1-E_2$, then $E=4L-3E_1-E_2-\cdots-E_9$
is an exceptional curve. Moreover, $p^*\Omega(1)|_C$ has the same splitting
as does $p^*\Omega(1)|_E$, since $C$ and $E$ both have the same
image $p(C)=p(E)$ in $\pr2$.

Thus an algorithm for computing the splitting for exceptional curves
handles all smooth rational $C$. 
\end{remark}

\setcounter{mythm}{0}
\subsection{The SHGH Conjecture}\label{SHGHsect}

Given any class $F\in \hbox{Cl}(X)$, there is a 
geometrically defined quantity $e(h^0,F)$ such that 
$h^0(X, F)\ge e(h^0,F)$ holds for general points $P_i$.
We now define this lower bound.

If $F\not\in\Psi_n$, then $h^0(X, F)=0$ by Lemma \ref*{Wlem},
and we set $e(h^0,F)=0$. If $F\in\Psi_n$, then we have the decomposition
$F=H+N$ given by Lemma \ref*{Wlem}(e), and we have
$h^0(X,F)=h^0(X,H)$. Since $H\cdot L\ge0$, we have
$h^0(X,H)\ge \hbox{max}(0, (H^2-K_X\cdot H)/2+1)$ and we set
$e(h^0,F)=\hbox{max}(0, (H^2-K_X\cdot H)/2+1)$.
Clearly, $h^0(X, F)\ge e(h^0,F)$ holds.
We can now state the SHGH Conjecture, which says
that equality in fact holds:

\addtocounter{mythm}{1}
\begin{conj}\label{SHGHconj}
We have  $h^0(X, F)= e(h^0,F)$,
where $F\in\hbox{Cl}(X)$ and $X$ is the blow up
of $\pr2$ at general points $P_i$.
\end{conj}

A version of the SHGH Conjecture for fat points follows from this.
Given a fat point subscheme $Z=m_1P_1+\cdots+m_nP_n$
supported at general points $P_i$, we define
$e(h_Z,t)$ to be $e(h^0, F_t(Z))$, where
$F_t(Z)=tL-m_1E_1-\cdots-m_nE_n$.

\addtocounter{mythm}{1}
\begin{conj}\label{SHGHconjfatpts}
We have  $h_Z(t)= e(h_Z,t)$,
where $Z=m_1P_1+\cdots+m_nP_n$
is a fat point scheme supported at general points $P_i$
of $\pr2$.
\end{conj}

To apply these conjectures, one must be able to compute
$e(h^0,F)$. We give two examples showing how to do so.

\addtocounter{mythm}{1}
\begin{example}\label{expex1}\rm
Suppose $F=tL-(77(E_1+\cdots+E_7)+44E_8+ 11E_9+ 11E_{10}+ 11E_{11})$.
To compute $e(h^0,F)$ for any given $t$, just mimic
the proof of Lemma \ref*{Wlem}(e). 
The idea is to find an element $w\in W$ such that
either $wF\cdot L$ is as small as possible, or 
$wF\cdot L<0$ or $wF\cdot (L-E_1)<0$.
For example, say $t=208$.
Apply $s_0$ to $F$ to get
$s_0F=185L-54E_1-54E_2-54E_3-77E_4-77E_5
-77E_6-77E_7-44E_8-11E_9-11E_{11}-11E_{11}$.
Permute the $E_i$ so that the coefficients
are nondecreasing, which gives
$F'=185L-77E_1-77E_2
-77E_3-77E_4-54E_5-54E_6-54E_7-44E_8-11E_9-11E_{11}-11E_{11}$.
This operation, taking $F$ to $F'$, is now repeated
until we obtain a class $F''$ such that either
$F''\cdot L<0$, or $F''\cdot (L-E_1)<0$,
or until $F''\cdot L\ge0$, $F''\cdot (L-E_1)\ge0$
and $F''\cdot r_0\ge0$. In this case
the class $F''$ we eventually end up with is
$-23L-8E_1+E_2+5E_3+5E_4+5E_5+5E_6+8E_7+8E_8+14E_9+17E_{10}+17E_{11}$,
hence $e(h^0,F)=0$, since $F''\cdot L<0$ (so $F''\not\in\Psi_{11}$).
If, for example, $t=209$, then the class we end up with is
$F''=11E_{11}$, so the decomposition of Lemma \ref*{Wlem}(d)
is $H=0$ and $N=F$, so $e(h^0,F)=1$. 
And if $t=210$, then we end up with $F''=
27L-8(E_1+\cdots+E_4)-5(E_5+\cdots+E_{11})$,
which is in $\Delta$, so $F=H$, $N=0$, and
$e(h^0,F)=(H^2-K_X\cdot H)/2+1$.
In this example, the SHGH Conjecture, that
$h^0(X, F)=e(h^0,F)$, in fact holds for all $t$.
It holds for $t<209$ since $F\not\in\Psi_n$
for those $t$. It holds for $t=209$,
since $wF=11E_{11}$ for some $w$, so $h^0(X, F)=h^0(X, 11E_{11})=1$.
And it holds for $t>209$ since for these cases
$F\in {\cal E}^*_{11}= W_{11} \Delta_{11}$, so
$F=H$, and $-K_X\cdot F\ge0$, so
$h^0(X, F)=(F^2-K_X\cdot F)/2+1$ by
Theorem 1.1 of \cite{refTrans} and semicontinuity of $h^0$.
(Macaulay 2 scripts for carrying out both the 
Lemma \ref*{Wlem}(d) decomposition and the Weyl group
calculations,
and a sample Macaulay 2 session demonstrating their use,
are included in Section A2.3 of the posted version
of this paper, \cite{refGHIpv}.)
\end{example}

\addtocounter{mythm}{1}
\begin{example}\label{expex2}\rm
Now consider $F=tL-50E_1-50E_2-38E_3-38E_4-26E_5-26E_6
-22E_7-18E_8-14E_9-14E_{10}$.
As in Example \ref*{expex1}, we have $e(h^0,F)=0$ for $t<102$, since
$F\not\in\Psi_{10}$. 
For $t=102$, we find a $w$ such that $wF=6L- 2(E_2+\cdots+E_8)+2E_9+6E_{10}$.
Thus the decomposition $F=H+N$ has $H=w^{-1}(6L- 2(E_2+\cdots+E_8))$ and
$N=w^{-1}(2E_9+6E_{10})$, where $w^{-1}$ can be performed 
by simply reversing the operations which gave $w$. What we find is
$H=38L-18E_1-18E_2-14E_3-14E_4-10E_5-10E_6-8E_7-8E_8-6E_9-6E_{10}$,
and $N=2C_1+6C_2$, where
$C_1=8L-4E_1-4E_2-3E_3-3E_4-2E_5-2E_6-E_7-2E_8-E_9-E_{10}$
and $C_2=8L-4E_1-4E_2-3E_3-3E_4-2E_5-2E_6-2E_7-E_8-E_9-E_{10}$.
Thus $h^0(X, F)=h^0(X, H)$, and it is known that $h^0(X,H)=4$.
For $t\ge103$, we have $F=H$ and $N=0$.
In fact, for $t=103$ we have $h^0(X,F)=92$, and for $t=104$ we have
$h^0(X,F)=197$.
(For the same reasons as in Example \ref*{expex1}, the SHGH Conjecture
holds for $F$ for all $t$.)
\end{example}

Given a fat point subscheme $Z\subset \pr2$ with general support, we can now
define the expected value $e(g_{{}_{\scriptscriptstyle\bullet}}(Z),i)$
of the Betti number $g_i(Z)$ for $i > \alpha+1$:

\addtocounter{mythm}{1}
\begin{defn}\label{expg}\rm
Let $F=F_{i-2}(Z)$; note that $F$ is effective, since
$i > \alpha+1$. Thus we have a decomposition
$F=H+N$, with $N=c_1C_1+\cdots+c_rC_r$, as in Lemma \ref*{Wlem}(e).
Let $m_i=\hbox{min}(c_i,L\cdot C_i)$. Let $M=m_1C_1+\cdots+m_rC_r$,
so $N-M$ is effective and $F=H+M+(N-M)$.
Then $\hbox{dim}\,\hbox{cok}\,\mu_{L+H+M}=\hbox{dim}\,\hbox{cok}\,\mu_{L+M}$
by Proposition \ref*{fitchettidea} assuming the SHGH Conjecture, so
$g_i(Z)=\hbox{dim}\,\hbox{cok}\,\mu_{L+F}=\hbox{dim}\,\hbox{cok}\,\mu_{L+M}
+ (h^0(X, 2L+F)-h^0(X, 2L+H+M))$ by Lemma \ref*{baselocuslem}, 
and $\hbox{dim}\,\hbox{cok}\,\mu_{L+M}\le
\sum_i {m_i-b_{C_i}\choose 2} + {m_i-a_{C_i}\choose 2}$ by 
Theorem \ref*{fitchettbounds} with equality assuming
Conjecture \ref*{boundseqconj}.
Thus we take $e(g_{{}_{\scriptscriptstyle\bullet}}(Z),i)$
to be $(h^0(X, 2L+F)-h^0(X,2L+H+M))+\sum_i {m_i-b_{C_i}\choose 2} + 
{m_i-a_{C_i}\choose 2}$. (In the notation of the proof of Theorem
\ref*{Bettibound}, $M$ is $N-N'$ and $m_i=c_i-c_i'$, 
so the upper bound in Theorem \ref*{Bettibound}
is by the proof of Theorem \ref*{Bettibound} exactly $e(g_{{}_{\scriptscriptstyle\bullet}}(Z),i)$.)
\end{defn}

The following result will be useful in our examples.
Versions of this result were proved in
\cite{refigp} and \cite{refFHH} and were the basis for the results
in \cite{refCJM} and \cite{refFHH}.

\addtocounter{mythm}{1}
\begin{lem}\label{qstarlstarlem}
Let $Z=m_1P_1+\cdots+m_nP_n$ be a fat point subscheme of
$\pr2$. Assume $F=F_k(Z)$ is the class of an effective divisor on $X$, and define
$h(F)=h^0(X,F)$, $l(F)=h^0(X, F-(L-E_1))$, $l^*(F)=h^1(X, F-(L-E_1))$,
$q(F)=h^0(X, F-E_1)$,
and $q^*(F)=h^1(X, F-E_1)$. Then
$$l(F)\le \hbox{dim ker }\mu_F \le l(F)+q(F).$$
If moreover
$h^1(X,F)=0$, then
$$k+2-2h(F)+l(F)\le \hbox{dim cok }\mu_F\le q^*(F)+l^*(F).$$
\end{lem}

Here are two examples showing explicitly how to compute expected values of the
graded Betti numbers.

\addtocounter{mythm}{1}
\begin{example}\label{boundsexample}\rm
Suppose we want to determine the graded Betti numbers
of the fat points subscheme $Z=m_1P_1+\cdots+m_nP_n\in \pr2$
where $n=10$ here, the points $P_i$ are general and the sequence of
multiplicities $m_i$ is $(50,50,38,38,26,26,22,18,14,14)$.
We found the Hilbert function in Example \ref*{expex2},
from which it follows that
$g_i=0$ except for $i=102$ (where we have $g_{102}=h_Z(102)=4$)
and possibly for $i=103$ and $i=104$.
We do not have an expected value for $g_{103}$, since
$103=\alpha(Z)+1$, but an ad hoc use of Lemma \ref*{qstarlstarlem}
gives $\hbox{dim}\,\hbox{ker}\,\mu_{F_{102}(Z)}\le l(F_{102}(Z))+q(F_{102}(Z))$,
where $l(F_{102}(Z))=h^0(X, F_{102}(Z)-(L-E_1))$ and $q(F_{102}(Z))=
h^0(X, F_{102}(Z)-E_1)$. Neither $F_{102}(Z)-(L-E_1)$ nor
$F_{102}(Z)-E_1$ is in $\Psi$, so $l=q=0$, so
$\mu_{F_{102}(Z)}$ is injective, hence
$g_{103}=\hbox{dim}\,\hbox{cok}\,\mu_{F_{102}(Z)}=h^0(X, F_{103}(Z))-3h^0(X,F_{102}(Z))
=92 - 3(4) = 80$. To compute $g_{104}$, recall in Example \ref*{expex2} we found
$F_{103}(Z)=H+2C_1+6C_2$. From
Lemma \ref*{splitlem} we have $a_{C_i}=b_{C_i}=4$ for $i=1,2$.
Since $h^1(X, H)=0$ holds here, by Proposition \ref*{fitchettidea}
we have $g_{104}=\hbox{dim}\,\hbox{cok}\,\mu_{F_{103}(Z)}=
\hbox{dim}\,\hbox{cok}\,\mu_{L+2C_1}+\hbox{dim}\,\hbox{cok}\,\mu_{L+6C_2}$.
By Theorem \ref*{fitchettbounds} we have
$\hbox{dim}\,\hbox{cok}\,\mu_{L+2C_1}=0$ and $\hbox{dim}\,\hbox{cok}\,\mu_{L+6C_2}=2$,
so $g_{104}=2$. For $i\ge104$, we have $h^1(X, F_i(Z))=0$, hence
$g_{i+1}=0$ by Corollary \ref*{CMlem}.

We can now write down a minimal free graded resolution for $I(Z)$.
It is $0\to M_1\to M_0\to I(Z)\to 0$, where
$M_0=R^{2}[-104]\oplus R^{80}[-103]\oplus R^4[-102]$, and 
from the Hilbert functions of $I(Z)$ and $M_0$ we now find 
$M_1=R^{16}[-105]\oplus R^{69}[-104]$.
\end{example}

\addtocounter{mythm}{1}
\begin{example}\label{rarecases}\rm
Consider $Z=48P_1+33P_2+33P_3+33P_4+32P_5+32P_6+32P_7+24P_8+16P_9$,
where the points $P_i$ are general.
Then $h_Z(t)=0$ for $t<98$, since $F_t(Z)\not\in\Psi$, and
$h_Z(t)={t+2\choose 2}-4879$ for $t\ge98$, since then $h_Z(t)>0$, $F_t(Z)\in {\cal E}^*$ and
the SHGH Conjecture is known to hold for $n=9$ general points (see \cite{refN2},
or use the fact that any nef divisor $F$ on a blow up of $\pr2$ at 9 general points 
has $-K_X\cdot F\ge0$ and apply the results of \cite{refTrans}). We also know that
$g_t=0$ for $t<98$, and $g_{98}=h_Z(98)=71$.
By Corollary \ref*{CMlem}, we have $g_t=0$ for $t>99$. In this case
$\alpha=98$, and we do not in general
have a conjectural value for $g_{\alpha+1}$, but in this case
there is an element $w\in W$ such that $F_{97}(Z)=L+H+N$
where $H=L-E_2-E_3-E_4$, $N=8E$, where $E$ is the exceptional curve
$12L-6E_1-4E_2-\cdots-4E_7-3E_8-2E_9$.
Since $h^1(X, H)=0$, reasoning as in the proof of
Proposition
\ref*{fitchettidea}, the cokernels of $\mu_{L+H+8E}$ and $\mu_{L+8E}$ are both
isomorphic to the cokernel of $\mu_{E;L+8E}$ and hence to each other.
By Lemma \ref*{splitlem} we have $a_e=6=b_E$, so by 
Theorem \ref*{fitchettbounds} we have $g_{99}=\hbox{dim}\,\hbox{cok}\,\mu_{L+8E}=2$.
The minimal free graded resolution for $I(Z)$ is thus
$0\to R^{28}[-100]\oplus R^{44}[-99]\to R^{2}[-99]\oplus R^{71}[-98]\to I(Z)\to0$.
\end{example}

\addtocounter{myapp}{1}
\setcounter{thm}{0}
\section{Computational Aspects}\label{algrthms}

In this section we discuss various computational aspects
of the problem of computing graded Betti numbers, partly
to explain how our geometric approach
can be used to make computer calculations more efficient,
and partly to give additional evidence in support
of Conjecture \ref*{boundseqconj}. 

\setcounter{mythm}{0}
\subsection{Splitting Types}\label{splittingtypes}

The first issue is the need to determine the splitting of the restriction
of $p^*\Omega_{\pr2}(1)$ to a smooth rational curve $C\subset X$,
where $p:X\to \pr2$ is a blow up of general points.
We are mainly interested in doing this for exceptional curves,
but it is of interest also to consider any smooth rational curve $C$.
However, doing so for exceptional curves suffices to do it
for all other smooth rational $C$ (see Remark \ref*{excsuffice}).

So suppose $E$ is an exceptional curve on $X$.
The simplest approach conceptually is to find a Cremona transformation $w$
of the plane that transforms $p(E)$ to a line $A$. Pick a general basis
$f_1,f_2,f_3$ of the linear forms on $\pr2$
and find their images $g_i=f_i\circ w^{-1}$ under $w$.
Given the equation of $A$, which is easy to get since $\hbox{deg}(A)=1$,
we can find the ideal $J'$ generated by the restrictions of the $g_i$ to $A$.
This ideal typically has base points; the ideal $J$ residual to the base points
is just the ideal generated by the restriction of the $f_i$ to $E$,
but regarding
$A=E$ as $\pr1$, we can find a minimal free resolution of the ideal
over $K[\pr1]$. The degree of the syzygy of least degree is $a_E$;
then $b_E=d-a_E$ where $d$ is the degree of $p(E)$.

It is not hard to convert this conceptual algorithm into code.
For speed, all of the actual symbolic operations should be done in
$K[\pr1]$ rather than in $K[\pr2]$. That is, one does not 
actually want to find the $g_i$ first and then restrict to the line.
A discussion for how to push all of the computation down to
$K[\pr1]$ is included in Section \ref*{algrthms}
of the posted version of the paper, \cite{refGHIpv}.
In addition, an explicit script that implements the computation
and which is very fast is included in Section A2.3
of the posted version.

It would be nice to be able to predict what $a_E$ and $b_E$
should be, based only on knowing $t$ and the $m_i$, given
a class $[E]=tL-m_1E_1-\cdots-m_nE_n$. This seems to be a difficult problem,
with Lemma \ref*{splitlem} being the main result. 
However, computational data suggests a possible new constraint
on the splitting type for a smooth rational curve $C$
with $C^2=1$. Using sequences analogous to $(\ddagger)$
we get homomorphisms
$$0 \to \hbox{cok }\mu_{1C} \to \hbox{cok } \mu_{2C} \to \hbox{cok }
\mu_{3C} \to\cdots \eqno(*).$$
The cokernel of $\hbox{cok }\mu_{iC} \to \hbox{cok } \mu_{(i+1)C}$
is $\hbox{cok }\mu_{C;(i+1)C}
\cong H^1(C, {\cal O}_C(i+1-a_C)\oplus{\cal O}_C(i+1-b_C))$.
For $r \ge b_C-2$, we thus get an upper
bound $\hbox{dim cok } \mu_{rC} \leq (a-1)(a-2)/2 + (b-1)(b-2)/2$ by adding up
the dimensions of the cokernels of $\hbox{cok } \mu_{C;(i+1)C}$ for $i\leq r$.

Moreover, $C=wL$ for some $w\in W(X)$, hence $C$ is
part of an exceptional configuration; i.e.,
$C = wL, C_1 = wE_1, \ldots, C_n = wE_n$.
When $r$ is big enough, $|rC-L|$ is
nonempty and the fixed part of $|rC - L|$ is $d_1C_1 + \cdots + d_nC_n$,
where $d_i = C_i\cdot L$. Thus, by Proposition \ref*{fitchettidea}, 
Conjecture \ref*{boundseqconj} 
and the SHGH Conjecture, we have $\hbox{dim cok }\mu_{rC} =
((a_1^2-a_1)/2 + (b_1^2-b_1)/2) + \cdots + ((a_n^2-a_n)/2 + (b_n^2-b_n)/2)$,
where $(a_i,b_i)$ is the splitting type for $C_i$.
This gives the inequality
$$((a_1^2-a_1)/2 + (b_1^2-b_1)/2) + \cdots + ((a_n^2-a_n)/2 + (b_n^2-b_n)/2)
\leq (a-1)(a-2)/2 + (b-1)(b-2)/2.\eqno(**)$$

This can be enough to determine $(a_C,b_C)$. For example,
let $C = 12L-5E_1-5E_2-5E_3-4E_4-4E_5-4E_6-4E_7-2E_8$.
So C comes from a plane curve of degree 12 with
three points of multiplicity 5, four of multiplicity 4 and one of
multiplicity 2.
By Lemma \ref*{splitlem}, $(a_C,b_C)$ is either $(6,6)$ or $(5,7)$.
Here are the $C_i$ and their types $(a_i,b_i)$.
For convenience we only give the coefficients of $L, -E_1,\ldots,-E_8$
followed by the splitting type:

$C_1 = (5; 2, 2, 2, 2, 1, 2, 2, 1), (2, 3)$

$C_2 = (5; 2, 2, 2, 2, 2, 1, 2, 1), (2, 3)$

$C_3 = (5; 2, 2, 2, 1, 2, 2, 2, 1), (2, 3)$

$C_4 = (5; 2, 2, 2, 2, 2, 2, 1, 1), (2, 3)$

$C_5 = (4; 2, 2, 2, 1, 1, 1, 1, 1), (2, 2)$

$C_6 = (3; 1, 1, 2, 1, 1, 1, 1, 0), (1, 2)$

$C_7 = (3; 2, 1, 1, 1, 1, 1, 1, 0), (1, 2)$

$C_8 = (3; 1, 2, 1, 1, 1, 1, 1, 0), (1, 2)$

Now $((a_1^2-a_1)/2 + (b_1^2-b_1)/2) + \cdots + ((a_n^2-a_n)/2 +
(b_n^2-b_n)/2)$ here is 21. But $(a-1)(a-2)/2 + (b-1)(b-2)/2 = 20$ if we use
$(a_C,b_C) = (6,6)$, so we see that $(a_C,b_C) = (5,7)$; this is also
what we get if we use Macaulay 2 \cite{refGS} to compute $(a_C,b_C)$,
using randomly chosen points (so this is a check but not a proof that $(a,b) = (5,7)$
here).

Also, it follows from $(*)$ and $(**)$ that
$\hbox{dim cok }\mu_{iC} = \hbox{dim cok } \mu_{C;C} +
\cdots + \hbox{dim cok } \mu_{C;iC}$
for all $i$ if $(**)$ is an equality.

Thus $(**)$ can, conjecturally, sometimes tell us both what the splitting type
is and what the cokernel is. Moreover, $(**)$ sometimes
also applies to exceptional
curves. In the example above, $E=C-E_9-E_{10}$ is exceptional and
has the same splitting type as does $C$,
so we get information about $E$ via
$C$.

We have applied our script for computing splitting types to
numerous examples.
In these examples, $(a_C,b_C)$ always made
$(a_C-1)(a_C-2)/2 + (b_C-1)(b_C-2)/2$ as small as possible subject to
$(**)$. (Moreover, in those cases where $(**)$ was not an equality,
it was off by exactly 1, and in those cases it always happened that $a_i=b_i$ for all $i$.)
This and Lemma \ref*{splitlem} lead us to make the following conjecture:

\addtocounter{mythm}{1}
\begin{conj}\label{splittingconj1}
Let $C=wL$ for some $w\in W(X)$, let $d=C\cdot L$,
let $m$ be the maximum
of $C\cdot E_1,\ldots,C\cdot E_n$ and let
$C_1=wE_1,\ldots, C_n=wE_n$. Then $(a_C,b_C)$
is the solution $(a,b)$ to $a\le b$,
$\hbox{min}(m,d-m)\le a\le d-m$ and $d=a+b$
which minimizes $(a-1)(a-2)/2 + (b-1)(b-2)/2$
subject to $(**)$.
\end{conj}

\setcounter{mythm}{0}
\subsection{Computational Evidence for Conjecture \ref*{boundseqconj}}\label{compevforconj}

There are 2051 exceptional classes of the form
$E=tL-m_1E_1-\cdots-m_nE_n$
with $1\le t\le 20$ (taking $n$ to be as large as necessary) and 
$m_1\ge\cdots\ge m_n\ge0$.
We have applied our splitting script
(using randomly chosen
points $P_i$ and working in characteristic 31991)
to  determine the splitting types of all 2051. It turned out that
Theorem \ref*{fitchettbounds}(b, c) implies for all but 25 of the 2051 cases
that  Conjecture \ref*{boundseqconj} holds for the given $E$.

Theorem \ref*{fitchettbounds}(b, c) does not apply
in the remaining 25 cases. Here we list these 25 cases,
giving $a_E$, $b_E$, $E\cdot L$ and $E\cdot E_i$ for all $i$ such that
$E\cdot E_i>0$:

\begin{verbatim}     
     5  8 13 5 5 5 5 5 5 4 1 1 1 1            8 11 19 8 8 7 7 7 6 6 3 2 1 1
     6  9 15 6 6 6 6 5 5 5 2 1 1 1            8 11 19 8 8 7 7 7 7 5 2 2 2 1
     6  9 15 6 6 6 6 6 5 4 1 1 1 1 1          8 11 19 8 8 8 7 6 6 6 2 2 2 1
     6 10 16 6 6 6 6 6 6 6 1 1 1 1 1          8 11 19 8 8 8 7 6 6 6 3 1 1 1 1
     7 10 17 7 7 6 6 6 6 6 2 2 2              8 11 19 8 8 8 7 7 6 5 2 2 1 1 1
     7 10 17 7 7 6 6 6 6 6 3 1 1 1            8 11 19 8 8 8 7 7 7 4 1 1 1 1 1 1 1
     7 10 17 7 7 7 6 6 6 5 2 2 1 1            8 11 19 8 8 8 8 6 6 5 2 1 1 1 1 1
     7 10 17 7 7 7 7 6 5 5 2 1 1 1 1          8 12 20 8 8 8 7 7 7 7 2 2 2 1
     7 10 17 7 7 7 7 6 6 4 1 1 1 1 1 1        8 12 20 8 8 8 7 7 7 7 3 1 1 1 1
     7 11 18 7 7 7 7 7 6 6 2 1 1 1 1          8 12 20 8 8 8 8 7 7 6 2 2 1 1 1
     7 11 18 7 7 7 7 7 7 5 1 1 1 1 1 1        8 12 20 8 8 8 8 8 6 6 2 1 1 1 1 1
     8 11 19 7 7 7 7 7 7 7 4 1 1 1            8 12 20 8 8 8 8 8 7 5 1 1 1 1 1 1 1
     8 11 19 8 7 7 7 7 7 6 3 2 2
\end{verbatim}

These cases can be checked by directly computing
$\hbox{dim}\,\hbox{cok}\,\mu_F$ for $F=L+iE$.
The critical value of $i$ is $i=b_E$ (if Conjecture \ref*{boundseqconj}
holds for $L+b_EE$,
it follows from the sequences $(\ddagger)$ that it holds for 
$L+iE$ for all $0\le i\le E\cdot L$). 
In principle, one can compute $\hbox{dim}\,\hbox{cok}\,\mu_F$
by computing a resolution of $I(Z)$ for
$Z=m_1P_1+\cdots+m_nP_n$ in the usual way, using Gr\"obner bases.
But for all but one of the 25 examples above we found this computation 
too large to successfully complete using this usual approach.
By taking a different geometrically inspired
approach we have in fact been able to check that 
Conjecture \ref*{boundseqconj} holds in these 24 other cases.
We now discuss this alternate approach
for computing the dimension of the cokernel of
$\mu_{L+mE}$, when $E$ is an exceptional curve and $0\le m\le L\cdot E$.
Since $\mu_{L+mE}$ and $\mu_{mE;L+mE}$ have isomorphic cokernels,
it is enough to compute the dimension of the cokernel of the latter, which 
it turns out one can do fairly efficiently, regarding the scheme structure of
$mE$ as Proj of a certain graded ring. So first we determine
this scheme structure.

Let $E\subset X$ be smooth and rational with $E^2=-1$.
First, we would like to know $\hbox{Pic}(mE)$, and to determine
$h^i(mE, {\cal F})$ for both $i=0$ and $i=1$ for every
line bundle ${\cal F}$ on $mE$.

Using exact sequences of the form
$0\to {\cal O}_X(dL+(j-m)E)\to {\cal O}_X(dL+jE)\to {\cal O}_{mE}(dL+jE)\to 0$
and
$0\to {\cal O}_X((d-1)L+jE)\to {\cal O}_X(dL+jE)\to {\cal O}_{L}(dL+jE)\to 0$
and the cohomology of divisors of the form $dL+jE$, which is known
for all $d$ and $j$, we find that: 
$h^0(mE, {\cal O}_{mE}(dL+jE)) - h^1(mE,{\cal O}_{mE}(dL+jE)) 
= {m+1\choose 2} + mt$ holds for all $t=E\cdot(dL+jE)$;
$h^0(mE, {\cal O}_{mE}(dL+jE))= {m+1\choose 2} + mt$ 
(and hence $h^1(mE, {\cal O}_{mE}(dL+jE))=0$)
and $h^0(mE, {\cal O}_{mE}(-dL-jE))= {m-t+1\choose 2}$ 
hold if $t\ge0$;
and it follows that $h^1(mE, {\cal O}_{mE}(-dL-jE))= 
{t\choose 2}$ holds if $0\le t\le m$,
and $h^1(mE, {\cal O}_{mE}(-dL-jE))=tm-{m+1\choose 2}$ holds 
if $t\ge m$.

In particular, $h^1(mE, {\cal O}_{mE})=0$, so it follows that
the inclusion $E\subset mE$ induces an isomorphism
$\hbox{Pic}(mE)\to\hbox{Pic}(E)={\bf Z}$ (see \cite{refArt},
(1.3) and (1.4)). It now follows for any line bundle $F$ on $X$ that 
$h^i(mE, {\cal O}_{mE}(F))$
depends only on $m$, $i$ and $E\cdot F$.
In particular, $h^i(mE, {\cal O}_{mE}(F))=h^i(mE, {\cal O}_{mE}(t))$,
where $t=E\cdot F$ so we define ${\cal O}_{mE}(1)$
to be ${\cal O}_{mE}(-E)$, and write $h^i(mE, t)$
for $h^i(mE, {\cal O}_{mE}(t))$.

To get the scheme structure, we can pick any $E$ which is 
convenient, since $mE$ is isomorphic for all $E$ on $X$. 
Let $\pi:S\to T$ be the blow up of the point
$P$ defined by $x=0=y$ in $T=\hbox{Spec}(K[x,y])$, and take 
$E=\pi^{-1}(P)$.
We can regard $S$ as $\hbox{Proj}_A(A[u,v]/(xv-uy))$, where
$A=K[x,y]$, and $E$ as $\hbox{Proj}_A(A[u,v]/((xv-uy)+(x,y)))$, 
which is just $E=\hbox{Proj}_K(K[u,v])$. 
Similarly, $mE$ is the fiber
$\pi^{-1}(mP)$, so $mE=\hbox{Proj}_A(A[u,v]/((xv-uy)+(x,y)^m))$.
Note that the ring $B=A[u,v]/((xv-uy)+(x,y)^m)$ is graded, with $x$
and $y$ of degree 0
and $u$ and $v$ of degree 1. The homogeneous component
$B_t$ for $t\ge 0$ is just the span of the monomials $x^iy^ju^rv^s$ with
nonnegative exponents where $r+s=t$, and $i+j<m$ (since we
have modded out by $(x,y)^m$). In fact $B_t$ can be identified
with $H^0(mE, t)$. However, $B$ has no components of negative degree,
whereas $h^0(mE, t)>0$ for $t>-m$. But note that $x/u=y/v$
on the open set of $mE$ where $u$ and $v$ are neither 0. Thus
taking $x/u$ where $u\ne0$ and $y/v$ where $v\ne0$ defines an element
$\epsilon\in H^0(mE, -1)$. Let $C=K[u,v,\epsilon]/(\epsilon^m)$.
This is graded if we take $\epsilon$ to have degree $-1$.
We have a graded injective ring homomorphism
$B\to C$ given by sending $x\mapsto \epsilon u$ and
$y\mapsto \epsilon v$, and now we can identify
$C_t$ with $H^0(mE, t)$ for all $t$. In particular,
the monomials of the form $\epsilon^iu^rv^s$
with nonnegative exponents satisfying $r+s=t+i$
and $i\le m-1$ give a basis for $H^0(mE, t)$,
and we can regard $mE$ as $\hbox{Proj}_K(C)$.

The map $\mu_{mE;L+mE}$ factors via the restriction 
$H^0(X, L) \to H^0(mE, {\cal O}_{mE}(L))$ through 
$$\mu:H^0(mE, {\cal O}_{mE}(L+mE))\otimes H^0(mE, {\cal O}_{mE}(L)) 
\to H^0(mE, {\cal O}_{mE}(2L+mE)).$$
Thus $\hbox{Im}\,\mu_{mE;L+mE}$ is the $K$-span of the product 
$(H^0(mE, {\cal O}_{mE}(L+mE)))(\hbox{Im}\,(H^0(X, L)\to H^0(mE, {\cal O}_{mE}(L))))$
in $H^0(mE, {\cal O}_{mE}(2L+mE))$.
It has the same dimension as does
$$(H^0(mE_1, {\cal O}_{mE_1}(L'+mE_1)))(\hbox{Im}\,(H^0(X, L')\to H^0(mE_1, {\cal O}_{mE_1}(L'))))$$
in $H^0(mE_1, {\cal O}_{mE_1}(2L'+mE_1))$, 
where $wL = L'$ for some appropriate Cremona 
transformation $w$ with $wE = E_1$.

Now choose coordinates on $\pr2$ such that $K[\pr2]=K[a,b,c]$, where 
$E_1$ is the blow up of the point $a=b=0$.
Let $d=E_1\cdot L'=E\cdot L$. Note that $L'-E_1$ is nef, so 
$0\le (L'-E_1)\cdot L=L'\cdot L-d$; thus $2d-1\le d-1+L\cdot L'$. Now, 
$H^0(mE_1, {\cal O}_{mE_1}(2L'+mE_1))=H^0(mE_1, {\cal O}_{mE_1}(2d-m))$,
so the monomials of the form $u^iv^j\epsilon^k$ with
$2d-m\le i+j\le 2d-1$ and $0\le k=i+j-(2d-m)$ give a basis.
Thus there is a surjective map of the homogeneous component
$(((a,b,c)^{2d-1})\cap((a,b)^{2d-m}))_{d-1+L\cdot L'}$
onto $H^0(mE_1, {\cal O}_{mE_1}(2d-m))$,
defined by sending $a^ib^jc^{d-1-i-j+L\cdot L'}$ to $u^iv^j\epsilon^{i+j-2d+m}$,
where $(a,b,c)$ denotes the ideal in $K[a,b,c]$ 
generated by $a$, $b$ and $c$. Moreover, the kernel
of this surjective map is spanned by the monomials
$a^ib^jc^{d-1-i-j+L\cdot L'}$ with $i+j\ge 2d$.

Note that the 
elements of $H^0(X, L')$, regarded as homogeneous 
polynomials in $K[a,b,c]$, are certain polynomials 
$f(a,b,c)$ of degree $L\cdot L'$ such that the terms 
of $f(a,b,1)$ of least degree have degree $d$.
The image of $f(a,b,c)$ in $H^0(mE_1, {\cal O}_{mE_1}(L'))$ 
(i.e., the restriction of $f$ to $mE_1$) is just what you 
get if you formally simplify $(\epsilon^d)(f(u/\epsilon^d,v/\epsilon^d,1))$.
Thus we have a surjection of $((a,b,c)^{d-1})\cap((a,b)^{d-m})$
onto $H^0(mE_1, {\cal O}_{mE_1}(L'+mE_1))$, 
defined by sending $a^ib^jc^{d-1-i-j}$ to $u^iv^j\epsilon^{i+j-d+m}$.
This gives a surjective map of the homogeneous component
$((((a,b,c)^{d-1})\cap((a,b)^{d-m}))H^0(X, L') + ((a,b)^{2d}))_{d-1+L\cdot L'}$
onto $\hbox{Im}\,\mu_{mE_1;L'+mE_1}$ with kernel $((a,b)^{2d})_{d-1+L\cdot L'}$. 
Thus $\hbox{dim}\,\hbox{Im}\,\mu_{mE_1;L'+mE_1}$ equals 
$\hbox{dim}\,((((a,b,c)^{d-1})\cap((a,b)^{d-m}))H^0(X, L') + ((a,b)^{2d}))_{d-1+L\cdot L'}
- \hbox{dim}\,((a,b)^{2d})_{d-1+L\cdot L'}$. 
The hardest part of this calculation is finding $H^0(X, L')$.
It can be done either using a Gr\"obner basis calculation
(but one much smaller than what is needed to calculate
the dimension of the image of $\mu_{mE;L+mE}$ directly)
or by applying $w$ to a basis for $H^0(X, L)$.

\end{document}